\documentclass{AIMS}
 \usepackage{amsmath,amssymb}
 \usepackage{mathrsfs}
 
\def\currentvolume{X}
 \def\currentissue{X}
  \def\currentyear{200X}
   \def\currentmonth{XX}
    \def\ppages{X--XX}

\newtheorem{theorem}{Theorem}
\newtheorem{corollary}{Corollary}

\newtheorem{lemma}{Lemma}
\newtheorem{proposition}{Proposition}

\theoremstyle{definition}

\newtheorem{remark}{Remark}

\newcommand{\eps}{\varepsilon}
\newcommand{\ep}{e^{-\frac{i\bc}{2}\int_{-\infty}^x\abs{\pn}^k dy}}

\newcommand{\norm}[1]{\left\|#1\right\|}
\newcommand{\abs}[1]{\left|#1\right|}

\newcommand{\Real}{\mathbb R}

\newcommand{\no}{\nonumber}

\newcommand{\Sz}{\mathscr{S}}
\newcommand{\F}{\mathscr{F}}

\newcommand{\ls}{\leqslant}
\newcommand{\gs}{\geqslant}

\newcommand\supp{\text{{\,\rm supp}\,}}

\newcommand{\bc}[1][\lambda]{#1}
\newcommand{\dx}{\langle D_x\rangle}
\newcommand{\pn}[1][\ll]{P_{#1 N}u}

\allowdisplaybreaks \numberwithin{equation}{section}
\numberwithin{theorem}{section} \numberwithin{lemma}{section}
\numberwithin{corollary}{section}
\numberwithin{proposition}{section} \numberwithin{remark}{section}
\title[Well-posedness for 1D DNLS]
      {Well-posedness for one-dimensional derivative nonlinear Schr\"odinger equations }

\author[C. C. Hao]{}


\email{hcc@amss.ac.cn}

\thanks{The author is supported in part by the Scientific Research
Startup Special Foundation on Excellent PhD Thesis and Presidential
Award of Chinese Academy of Sciences, NSFC (Grant No. 10601061), and
the Innovation Funds of AMSS, CAS of China.}

\begin{document}
\maketitle

\centerline{\scshape Chengchun Hao}
\medskip
{\footnotesize
 \centerline{Institute of Mathematics}
  \centerline{Academy of Mathematics \& Systems Science, CAS}
  \centerline{Beijing 100080, P. R. China}
} 

\medskip

 \centerline{(Communicated by Gigliola Staffilani)}
 \medskip

\begin{abstract}
In this paper, we investigate the one-dimensional derivative
nonlinear Schr\"odinger equations of the form
$iu_t-u_{xx}+i\lambda\abs{u}^k u_x=0$ with non-zero $\lambda\in
\Real$ and any real number $k\gs 5$. We establish the local
well-posedness of the Cauchy problem with any initial data in
$H^{1/2}$ by using the gauge transformation and the Littlewood-Paley
decomposition.
\end{abstract}

\section{Introduction}

In the present paper, we consider the following Cauchy problem for
the derivative nonlinear Schr\"odinger equation
\begin{align}\label{dnls}
    &iu_t-u_{xx}+i\bc\abs{u}^k u_x=0,\quad (t,x)\in\Real^2,\\
    &u(0,x)=u_0(x),\label{data}
\end{align}
where $u=u(t,x): \Real^2\to \mathbb{C}$ is a complex-valued wave
function, both $\bc\neq 0$ and $k\gs 5$ are  real numbers.

A great deal of interesting research has been devoted to the
mathematical analysis for the derivative nonlinear Schr\"odinger
equations
\cite{CKSTT01,CKSTT02,FT80,FT81,Gr00,Ha93,HO92,HO94,KPV93,Oz96,Ta99}.
In \cite{KPV93}, C. E. Kenig, G. Ponce and L. Vega studied the local
existence theory for the Cauchy problem of the derivative nonlinear
Schr\"odinger equations
\begin{align*}
    iu_t+u_{xx}+f(u,\bar{u},u_x, \bar{u}_x)=0, \quad (t,x)\in \Real^2,
\end{align*}
 with small data $u(0,x)=u_0(x)$ in $H^{7/2}$ where $f$ is a
polynomial having no constant or linear terms with the lowest order
term of degree being greater than or equal to $3$. Subsequently, it
was improved to $H^3$ by N. Hayashi and T. Ozawa \cite{HO94}.

If the nonlinearity consists mostly of the conjugate wave $\bar{u}$,
then it can be done much better.  In the case $f=(\bar{u}_x)^k$, A.
Gr\"uenrock, in \cite{Gr00}, obtained local well-posedness when
$s>s_c=3/2-1/(k-1)$, $s\gs 1$, and $k\gs 2$ was an integer. In
particular, the global well-posedness in $H^1$ is obtained when
$f=i(\bar{u}_x)^2$ with the help of the Bourgain spaces (cf.
\cite{Bou93,Taobook}).

In \cite{Ta99}, H. Takaoka discussed the derivative nonlinear
Schr\"odinger equation of the form
\begin{align*}
    u_t-iu_{xx}+|u|^2u_x=0,\quad (t,x)\in\Real^2,
\end{align*}
and obtained the local well-posedness in $H^s$ for $s\gs 1/2$ by
performing a fixed point argument in an adapted Bourgain space
$X_{s,b}$ which yields a $C^\infty$-solution map.

A very similar equation to \eqref{dnls} is the generalized
Benjamin-Ono equation
\begin{align}\label{GBO}
    u_t+\mathcal {H}u_{xx}\pm u^ku_x=0, \quad (t,x)\in\Real^2,
\end{align}
where $u$ is a real-valued function, $\mathcal {H}$ is the Hilbert
transformation defined by
\begin{align*}
    \mathcal{H}f(x)=-i\int_\Real
    e^{ix\xi}\text{sgn}(\xi)\hat{f}(\xi)d\xi,
\end{align*}
and $k\gs 2$ is an integer, the symbol $\hat{\cdot}$ (or $\F$)
denotes the spatial Fourier transform. For this equation, L. Molinet
and F. Ribaud \cite{MR04a,MR04b} obtained the local well-posedness
in the Sobolev space $H^s$ for $s>1/2$ if $k=2,4$, $s\gs 3/4$ if
$k=3$ and $s\gs 1/2$ if $k\gs 5$ by using Tao's gauge
transformation. In \cite{KT05}, C. E. Kenig and H. Takaoka have
shown the global well-posedness for the case $k=2$ in $H^s$ for
$s\gs 1/2$ by combining the gauge transformation with a
Littlewood-Paley decomposition and following the compactness
argument with  a priori estimates with the help of the preservation
of the Hamiltonian and the $L^2$-mass.

In the present paper, we shall generalize the above results to the
derivative nonlinear Schr\"odinger equation with $k\gs 5$ by using
some ideas in \cite{KT05}. However, we have to reconstruct new and
complicated estimates for the case $k\gs 5$ which is quite different
from the case $k=2$.

We first state the main result of this paper as follows, though we
shall define later the function space $X_T$ at the end of this
section.

\begin{theorem}\label{thm}
For any $u_0\in H^{1/2}$, there exist a $T=T(\norm{u_0}_{H^{1/2}})$
and a unique solution $u$ of \eqref{dnls}-\eqref{data} satisfying
\begin{align*}
    u\in C([-T,T];H^{1/2})\cap X_T.
\end{align*}
\end{theorem}

For convenience, we now introduce some notations. For nonnegative
real numbers $A$, $B$, we use $A\lesssim B$ to denote $A\ls CB$ for
some $C>0$ which is independent of $A$ and $B$. $A\sim B$ means
$A\lesssim B\lesssim A$, and $A\ll B$ denotes $A\ls CB$ for some
small $C>0$ which is also independent of $A$ and $B$.

To give the Littlewood-Paley decomposition, let $\psi$ be a fixed
even $C^\infty$ function with a compact support,
$\supp\psi\subset\{\abs{\xi}<2\}$, and $\psi(\xi)=1$ for
$\abs{\xi}\ls 1$. Define $\varphi(\xi)=\psi(\xi)-\psi(2\xi)$. Let
$N$ be a dyadic number of the form $N=2^j$, $j\in
\mathbb{N}\cup\{0\}$ or $N=0$. Writing
$\varphi_N(\xi)=\varphi(\xi/N)$ for $N\gs 1$, we define the
convolution operator $P_N$ by $P_Nu=\check{\varphi}_N*u$, where the
symbol $\check{\cdot}$ (or $\F^{-1}$) denotes the spatial Fourier
inverse transform. We define the function $\varphi_0$ by
$\varphi_0(\xi)=1-\sum_N \varphi_N(\xi)$ and denote
$P_0u=\check{\varphi}_0*u$. Then we introduce a spatial
Littlewood-Paley decomposition \cite{Stein}
\begin{align*}
    \sum_N P_N=I.
\end{align*}
Throughout this paper, we often use the Littlewood-Paley theorem
(cf. \cite{Stein,Taobook})
\begin{align*}
    \norm{\left(\sum_N\abs{P_N\phi}^2\right)^{1/2}}_{L^p}\sim
    \norm{\phi}_{L^p},
\end{align*}
for $1<p<\infty$. We also use more general operators $P_{\ll N}$ and
 $P_{\lesssim N}$ which are defined by
\begin{align*}
    P_{\ll N}=\sum_{M\ll N}P_M,\quad P_{\lesssim N}=\sum_{M\lesssim
    N}P_M,
\end{align*}
and $P_{\gg N}$, $P_{\gtrsim N}$ and $P_{\sim N}$ which can be
defined in a similar way. The Littlewood-Paley operators commute
with derivative operators (including $\abs{\nabla}^s$ and
$i\partial_t-\partial_{xx}$), the propagator
$S(t)=e^{-it\partial_x^2}$, and conjugation operations, are
self-adjoint, and are bounded on every Lebesgue space $L^p$ and
homogeneous Sobolev space $\dot{H}^s$ if $1\ls p\ls \infty$.
Furthermore, they obey the following Sobolev and Bernstein estimates
for $\Real$ with $s\gs 0$ and $1\ls p\ls \infty$ (which is similar
to those of three dimensions \cite{CKSTT04}):
\begin{align*}
    \norm{P_{\gs N}f}_{L^p}\lesssim &N^{-s}\norm{\abs{\nabla}^s P_{\gs
    N}f}_{L^p},\\
    \norm{P_{\ls N}\abs{\nabla}^s f}_{L^p}\lesssim & N^s\norm{P_{\ls
    N}f}_{L^p},\\
    \norm{P_N\abs{\nabla}^{\pm s}f}_{L^p}\lesssim &N^{\pm
    s}\norm{P_Nf}_{L^p},
\end{align*}
which can be verified by combining the Bernstein multiplier theorem
\cite{BL} and the interpolation theorem of Sobolev spaces.

We define the Lebesgue spaces $L_T^qL_x^p$ and $L_x^p L_T^q$ by the
norms
\begin{align*}
    \norm{f}_{L_T^qL_x^p}=\norm{\norm{f}_{L_x^p(\Real)}}_{L_t^q([0,T])},
    \quad
    \norm{f}_{L_x^pL_T^q}=\norm{\norm{f}_{L_t^q([0,T])}}_{L_x^p(\Real)}.
\end{align*}
In particular, we abbreviate $L_T^qL_x^p$ or $L_x^p L_T^q$ as
$L_{x,T}^p$ in the case $p=q$.

We also use the elementary inequality \cite{CKSTT04}
\begin{align*}
    \norm{\left(\sum_{N}\abs{f_N}^2\right)^{1/2}}_{L_T^qL_x^p}\ls
    \left(\sum_N \norm{f_N}_{L_T^qL_x^p}^2\right)^{1/2},
\end{align*}
for all $2\ls q,p\ls \infty$ and arbitrary functions $f_N$, and the
dual version
\begin{align*}
    \left(\sum_N \norm{f_N}_{L_T^{q'}L_x^{p'}}^2\right)^{1/2}
    \ls\norm{\left(\sum_{N}\abs{f_N}^2\right)^{1/2}}_{L_T^{q'}L_x^{p'}},
\end{align*}
where $p'$ is the conjugate number of $p$ given by $1/p+1/p'=1$. It
is easy to verify that they also hold if we replace the norm
$L_T^qL_x^p$ by the norm $L_x^pL_T^q$ in both side of the above
inequalities.

Let $\langle\cdot\rangle=(1+\abs{\cdot}^2)^{1/2}$. We use the
fractional differential operators $D_x^s$ and $\langle D_x\rangle^s$
defined by
\begin{align*}
    D_x^s f=\F^{-1}\abs{\xi}^s\F f,
    \quad \langle D_x\rangle^s f=\F^{-1}\langle \xi\rangle^s\F f.
\end{align*}
Thus, we can introduce the resolution space. For $T>0$,  we define
the function space $X_T$ in a similar way as in \cite{KT05} by
\begin{align*}
    X_T:=\{u\in\Sz'((-T,T)\times\Real) : \norm{u}_{X_T}<\infty\},
\end{align*}
where
\begin{align*}
    \norm{u}_{X_T}=&\norm{u}_{L_T^\infty H^{1/2}}
    +\left(\sum_{N}\norm{\partial_x P_N u}_{L_x^\infty
    L_T^2}^2\right)^{1/2}\\
    &+\left(\sum_{N}\norm{P_N u}_{L_x^2
    L_T^\infty}^2\right)^{1/2}+\left(\sum_{N}\norm{\langle D_x\rangle^{\frac{1}{4}} P_N
    u}_{L_x^4
    L_T^\infty}^2\right)^{1/2}.
\end{align*}

\section{Gauge transformation}

We transform the equation \eqref{dnls} by introducing the following
complex-valued function $v_N: \Real^2\to \mathbb{C}$ for a dyadic
number $N$ given by
\begin{align}\label{transf}
    v_N(t,x)=e^{-\frac{i\bc}{2}\int_{-\infty}^x\abs{\pn(t,y)}^k dy}\pn[].
\end{align}
By computation, we have
\begin{align}\label{tran1}
    i\partial_t v_N-\partial_x^2
    v_N=&-i\bc\ep\left[P_N(\abs{u}^ku_x)-\abs{\pn}^k\pn[]_x\right]\no\\
    & -\frac{i\bc}{2}
    \ep\pn[](i\partial_t-\partial_x^2)\int_{-\infty}^x\abs{\pn(t,y)}^k
    dy\no\\
    &+\frac{\bc^2}{4}\ep\abs{\pn}^{2k}\pn[].
\end{align}

For the second term, we integrate by parts and have
\begin{align*}
    &(i\partial_t-\partial_x^2)\int_{-\infty}^x\abs{\pn(t,y)}^k
    dy\\
    =&i\int_{-\infty}^x
    \frac{k}{2}\abs{\pn}^{k-2}\left(\partial_t\pn\overline{\pn}+\pn\partial_t\overline{\pn}\right)dy\\
    &-\frac{k}{2}\abs{\pn}^{k-2}\left(\partial_x\pn\overline{\pn}+\pn\partial_x\overline{\pn}\right)\\
    =&\int_{-\infty}^x
    \frac{k}{2}\abs{\pn}^{k-2}\left(\overline{\pn}P_{\ll
    N}(u_{xx}-i\bc\abs{u}^ku_x)\right.\\
    &\left.\qquad\qquad\qquad\qquad\qquad -\pn P_{\ll N}(\bar{u}_{xx}+i\bc\abs{u}^k\bar{u}_x)\right)dy\\
    &-\frac{k}{2}\abs{\pn}^{k-2}\left(\pn_x\overline{\pn}+\pn\overline{\pn_x}\right)\\
    =&\int_{-\infty}^x \frac{k(k-2)}{4}
    \abs{\pn}^{k-4}\left[(\overline{\pn_y}\pn)^2-(\pn_y\overline{\pn})^2\right]dy\\
    &-\int_{-\infty}^x\frac{i\bc k}{2}\abs{\pn}^{k-2}P_{\ll N}\abs{u}^k\left(u_x+\bar{u}_x
    \right)dy
    -k\abs{\pn}^{k-2}\pn\overline{\pn_x}
\end{align*}
Thus, $v_N$ obeys the following differential-integral equation
\begin{align}\label{dnlsnew}
    &i\partial_t v_N-\partial_x^2 v_N(t,x)\no\\
    =&-i\bc\ep\left[P_N(\abs{u}^ku_x)-\abs{\pn}^k\pn[]_x\right]\no\\
    & -\frac{i\bc k(k-2)}{8}\ep\pn[]\int_{-\infty}^x
    \abs{\pn}^{k-4}\cdot\no\\
    &\qquad\qquad\left[(\overline{\pn_x}\pn)^2-(\pn_x\overline{\pn})^2\right]dy\no\\
    &-\frac{\bc^2 k}{4}\ep\pn[]\int_{-\infty}^x\abs{\pn}^{k-2}P_{\ll N}
    \abs{u}^k\left(u_x+\bar{u}_x \right)dy\no\\
    &+\frac{i\bc k}{2}\ep\abs{\pn}^{k-2}\pn[]\pn\overline{\pn_x}\no\\
    &+\frac{\bc^2}{4}\ep\abs{\pn}^{2k}\pn[]\no\\
    \equiv &
    I_{N,1}(t,x)+I_{N,2}(t,x)+I_{N,3}(t,x)+I_{N,4}(t,x)+I_{N,5}(t,x).
\end{align}
The equivalent integral equation reads
\begin{align}\label{inte}
    v_N(t)=&S(t)e^{-\frac{i\bc}{2}\int_{-\infty}^x\abs{\pn_0(y)}^k
    dy}\pn[]_0\no\\
    &-i\int_0^t S(t-\tau)[I_{N,1}+I_{N,2}+I_{N,3}+I_{N,4}+I_{N,5}](\tau)
    d\tau.
\end{align}

\section{Preliminaries}

In order to prove the a priori estimate for the equation of $v_N$,
we need the linear estimates associated with the one-dimensional
Schr\"odinger equation. We first recall the Strichartz estimates,
smoothing effects and maximal function estimates. For the proofs,
one can see \cite{KPV93,KT05}.

\begin{lemma}\label{lem4}
For all $\phi\in\Sz(\Real)$, $\theta\in [0,1]$ and $T\in(0,1)$,
\begin{align}\label{eq1}
    \norm{S(t)\phi}_{L_T^{\frac{4}{\theta}}L_x^{\frac{2}{1-\theta}}}\lesssim
    &\norm{\phi}_{L^2},\\
    \norm{S(t)P_N\phi}_{L_x^{\frac{2}{1-\theta}}
    L_T^{\frac{2}{\theta}}}\lesssim & \langle
    N\rangle^{\frac{1}{2}-\theta}\norm{\phi}_{L^2},\label{eq2}\\
    \norm{S(t)\phi}_{L_x^4 L_T^\infty}\lesssim &
    \norm{\phi}_{\dot{H}^{\frac{1}{4}}}.\label{eq3}
\end{align}
\end{lemma}

We also need the $L_T^qL_x^p$ and $L_x^p L_T^q$ estimates for the
linear operator $f\mapsto \int_0^t S(t-\tau) f(\tau)d\tau$. For the
proofs, one can see \cite{KT05}.

\begin{lemma}\label{lem5}
For $f\in\Sz(\Real^2)$, $\theta\in [0,1]$ and $T\in (0,1)$,
\begin{align}\label{eq4}
    \norm{\int_0^t
    S(t-\tau)f(\tau)d\tau}_{L_T^{\frac{4}{\theta}}L_x^{\frac{2}{1-\theta}}}\lesssim
    &\norm{f}_{L_T^{\left(\frac{4}{\theta}\right)'}L_x^{\left(\frac{2}{1-\theta}\right)'}},\\
    \norm{\dx^{\frac{\theta}{2}}\int_0^t S(t-\tau) f(\tau)d\tau}_{L_T^\infty
    L_x^2}\lesssim & \norm{f}_{L_x^{p(\theta)}L_T^{q(\theta)}},\label{eq5}\\
    \norm{\dx^{\frac{\theta}{2}}\int_0^t S(t-\tau)P_N f(\tau)d\tau}_{L_x^2L_T^\infty
    }\lesssim & \langle
    N\rangle^{\frac{1}{2}}\norm{f}_{L_x^{p(\theta)}L_T^{q(\theta)}},\label{eq7}\\
    \norm{\dx^{\frac{\theta}{2}-\frac{1}{4}}\int_0^t S(t-\tau)
    f(\tau)d\tau}_{L_x^4 L_T^\infty}\lesssim
    &\norm{f}_{L_x^{p(\theta)}L_T^{q(\theta)}},\label{eq8}\\
    \norm{\dx^{\frac{1}{2}}\int_0^t S(t-\tau)f(\tau)d\tau}_{L_x^\infty
    L_T^2}\lesssim & \norm{f}_{L_T^{1}L_x^{2}},\label{eq9}\\
    \norm{\int_0^t S(t-\tau)P_N
    f(\tau)d\tau}_{L_x^{\frac{2}{\theta}}L_T^{\frac{2}{1-\theta}}}\lesssim
    & \langle N\rangle^{\frac{1}{2}-\theta}\norm{f}_{L_T^1L_x^2},\label{eq10}\\
   \norm{\int_0^t S(t-\tau)f(\tau)d\tau}_{L_x^4 L_T^\infty}\lesssim
    & \norm{f}_{L_T^1\dot{H}_x^{\frac{1}{4}}},\label{eq11}\\
    \norm{\partial_x \int_0^t S(t-\tau)f(\tau)d\tau}_{L_x^\infty
    L_T^2}\lesssim & \norm{f}_{L_x^1 L_T^2},\label{eq12}
\end{align}
where $p'$ is the conjugate number of $p\in[1,\infty]$, i.e.
$1/p+1/p'=1$, and
\begin{align*}
    \frac{1}{p(\theta)}=\frac{3+\theta}{4}, \quad
    \frac{1}{q(\theta)}=\frac{3-\theta}{4}.
\end{align*}

\end{lemma}

Next, we recall the Leibniz' rule for a product of the form
$e^{iF}g$ where $F$ is the spatial primitive of some function $f$.
For the proof, we refer to \cite{KT05,MR04b}.

\begin{lemma}[{\cite[Lemma 3.5]{KT05}}]\label{lem1}
Let $\alpha\in (0,1)$, $p$, $p_1$, $p_2$, $q$, $q_1\in(1,\infty)$,
$q_2\in(0,\infty]$ with $\frac{1}{p}=\frac{1}{p_1}+\frac{1}{p_2}$,
$\frac{1}{q}=\frac{1}{q_1}+\frac{1}{q_2}$, and let
$F(t,x)=\int_{-\infty}^x f(t,y) dy$, with real-valued function $f$.
Then
\begin{align*}
    \norm{D_x^\alpha(e^{iF}g)}_{L_x^p L_T^q}\lesssim \norm{f}_{L_x^{p_1}
    L_T^{q_1}}\norm{g}_{L_x^{p_2}L_T^{q_2}}+\norm{\langle D_x\rangle^\alpha g}_{L_x^p
    L_T^q}.
\end{align*}
\end{lemma}

\section{Bilinear estimates}

In this section, we prove the following space-time estimate which is
crucial to the proof of the nonlinear estimates.

\begin{proposition}\label{prop1}
Let $u\in H^\infty$ and $p\gs 4$ be a real number. Then we have
\begin{align}\label{eqprop1}
    \norm{u\bar{u}_x}_{L_x^p L_T^2}\lesssim T^{\frac{1}{2}}\norm{u}_{X_T}^2
    +(1+T^{\frac{1}{4}}\norm{u}_{X_T})\norm{u}_{X_T}\norm{P_{\gg
    1}u}_{X_T}.
\end{align}
\end{proposition}

\begin{proof}
By the Littlewood-Paley decomposition, we can write
\begin{align}\label{eqprop11}
    \norm{u\bar{u}_x}_{L_x^p L_T^2}=&\norm{\sum_{N_1,N_2}P_{N_1}uP_{N_2}\bar{u}_x}_{L_x^p L_T^2}\no\\
    \lesssim &\norm{\sum_{N_1\sim N_2}P_{N_1}uP_{N_2}\bar{u}_x}_{L_x^p
    L_T^2}
    +\norm{\sum_{N_1\ll N_2}P_{N_1}uP_{N_2}\bar{u}_x}_{L_x^p
    L_T^2}\\
    &+\norm{\sum_{N_1\gg N_2}P_{N_1}uP_{N_2}\bar{u}_x}_{L_x^p L_T^2}\no\\
    =:&I_1+I_2+I_3.
\end{align}

Now, we derive the estimates for $I_1$, $I_2$ and $I_3$,
respectively.

From the H\"older inequality, the Bernstein type inequalities and
the real interpolation theorem, we have
\begin{align*}
    I_1\lesssim &\sum_{N_1\sim N_2}\norm{P_{N_1}u}_{L_x^{2p}
    L_T^{4}} \norm{P_{N_2}u_x}_{L_x^{2p}
    L_T^4}
    \lesssim \sum_{N_1\sim N_2}\norm{P_{N_1}u}_{L_x^{2p}
    L_T^{4}} N_2\norm{P_{N_2}u}_{L_x^{2p}
    L_T^4}\\
    \lesssim &\sum_{N_1\sim N_2}\norm{D_x^{1/2} P_{N_1}u}_{L_x^{2p}
    L_T^{4}} \norm{D_x^{1/2} P_{N_2}u}_{L_x^{2p}
    L_T^4}
    \lesssim \sum_{N}\norm{D_x^{1/2} P_{N}u}_{L_x^{2p}
    L_T^{4}}^2\\
    \lesssim &\sum_{N}\norm{D_x^{1/2} P_{N}(P_{\lesssim 1}u+P_{\gg 1}u)}_{L_x^{2p}
    L_T^{4}}^2\\
    \lesssim &\sum_{N}\norm{D_x^{1/2} P_{N}P_{\lesssim 1}u}_{L_x^{2p}
    L_T^{4}}^2+\sum_{N}\norm{D_x^{1/2} P_{N}P_{\gg 1}u}_{L_x^{2p}
    L_T^{4}}^2\\
    \lesssim &\sum_{N}\norm{P_{N}P_{\lesssim 1}u}_{L_x^{2p}
    L_T^{4}}^2+\sum_{N}N\norm{ P_{N}P_{\gg 1}u}_{L_x^{2p}
    L_T^{4}}^2\\
    \lesssim &\sum_{N\lesssim 1}\norm{P_{N}u}_{L_x^{2p}
    L_T^{4}}^2+\sum_{N}N\norm{ P_{N}P_{\gg 1}u}_{L_x^{p}
    L_T^{\infty}}\norm{ P_{N}P_{\gg 1}u}_{L_x^{\infty}
    L_T^{2}}.\\
\intertext{Applying the Sobolev embedding theorem and the H\"older
inequality to the first term, and Bernstein estimates to the second
term, we can see that it is bounded by}
    \lesssim &T^{1/2}\norm{u}_{X_T}^2+\sum_{N}\norm{ P_{N}P_{\gg 1}u}_{L_x^{p}
    L_T^{\infty}}\norm{\partial_xP_{N}P_{\gg 1}u}_{L_x^{\infty}
    L_T^{2}}.\\
\intertext{By the Cauchy-Schwartz inequality and the Sobolev
embedding theorem(i.e. $H_4^{1/4} \subset H_4^{1/4-1/p}\subset L^p$
for the real number $p\gs 4$), we can bound it by}
    \lesssim &T^{1/2}\norm{u}_{X_T}^2+\norm{P_{\gg
    1}u}_{X_T}^2.
\end{align*}

For $I_2$ or $I_3$, it is suffice to consider one of them, e.g.
$I_2$, in view of symmetry. For $N_1\ll N_2$, we have
\begin{align*}
    P_{N_1}uP_{N_2}\bar{u}_x=\tilde{P}_{N_2}(P_{N_1}uP_{N_2}\bar{u}_x),
\end{align*}
where $\tilde{P}_{N}=\sum_{j=-2}^2 P_{2^j N}$. We split these into
three cases, i.e. $N_1\lesssim 1\ll N_2$, $N_1\ll N_2\lesssim 1$ and
$1\ll N_1\ll N_2$. For the case $N_1\lesssim 1\ll N_2$, from the
H\"older inequality and the Littlewood-Paley theorem, we can get
\begin{align}\label{eqprop12}
    &\norm{\sum_{N_1\lesssim 1\ll N_2}P_{N_1}uP_{N_2}\bar{u}_x}_{L_x^p
    L_T^2}\no\\
    = &\norm{P_{\lesssim 1}uP_{N_2\gg 1}\bar{u}_x}_{L_x^p L_T^2}
    \lesssim \norm{P_{\lesssim 1}u}_{L_x^p L_T^\infty}\norm{P_{\gg 1}u_x}_{L_x^\infty L_T^2}\no\\
    \lesssim & \norm{\left(\sum_{M}\abs{P_MP_{\lesssim 1}u}^2\right)^{1/2}}_{L_x^p L_T^\infty}
    \norm{\left(\sum_{M}\abs{P_MP_{\gg 1}u_x}^2\right)^{1/2}}_{L_x^\infty L_T^2}\no\\
    \lesssim & \left(\sum_{M\lesssim 1}\norm{P_Mu}_{L_x^p L_T^\infty}^2\right)^{1/2}
    \left(\sum_{M}\norm{P_MP_{\gg 1}u_x}_{L_x^\infty L_T^2}^2\right)^{1/2}\no\\
    \lesssim &\norm{u}_{X_T}\norm{P_{\gs 1}u}_{X_T}.
\end{align}
For the case $N_1\ll N_2\lesssim 1$, we have, by the H\"older
inequality and the Littlewood-Paley theorem, that
\begin{align}
    &\norm{\sum_{N_1\ll N_2\lesssim 1}P_{N_1}uP_{N_2}\bar{u}_x}_{L_x^p L_T^2}
    =  \norm{\sum_{N_2\lesssim 1}P_{\ll N_2}uP_{N_2}\bar{u}_x}_{L_x^p L_T^2}\no\\
    \lesssim &T^{1/2}\norm{\left(\sum_{N_2\lesssim 1}\abs{P_{\ll
    N_2}u}^2\right)^{1/2}}_{L_x^{2p} L_T^\infty}
    \norm{\left(\sum_{N_2 \lesssim 1}\abs{P_{N_2}\bar{u}_x}^2\right)^{1/2}}_{L_x^{2p} L_T^\infty}\no \\
    \lesssim &T^{1/2}\left(\sum_{N_2\lesssim 1}\norm{P_{\ll
    N_2}u}_{L_x^{2p}
    L_T^\infty}^2\right)^{1/2}\norm{u}_{X_T}.\label{eqprop15}
\end{align}
For $N_2\lesssim 1$, we have, by the Sobolev embedding theorem, that
\begin{align*}
    &\norm{P_{\ll N_2}u}_{L_x^{2p} L_T^\infty}
    \sim  \norm{\left(\sum_{M}\abs{P_MP_{\ll N_2}u}^2\right)^{1/2}}_{L_x^{2p} L_T^\infty}
    \lesssim \left(\sum_{M\ll N_2}\norm{P_Mu}_{L_x^{2p}
    L_T^\infty}^2\right)^{1/2}\\
    \lesssim &\left(\sum_{M\ll N_2}N_2^{2\eps}M^{-2\eps}\norm{P_Mu}_{L_x^{2p}
    L_T^\infty}^2\right)^{1/2}
    \lesssim \left(\sum_{M\ll N_2}N_2^{2\eps}\norm{D_x^{-\eps}P_Mu}_{L_x^{2p}
    L_T^\infty}^2\right)^{1/2}\\
    \lesssim &N_2^\eps \left(\sum_{M\lesssim 1}\norm{P_Mu}_{L_x^{2}
    L_T^\infty}^2\right)^{1/2}\lesssim N_2^\eps\norm{u}_{X_T},
\end{align*}
where $\eps=(p-1)/2p$. Thus, \eqref{eqprop15} can be bounded by
\begin{align*}
    \lesssim T^{1/2}\left(\sum_{N_2\lesssim
    1}N_2^{2\eps}\right)^{1/2}\norm{u}_{X_T}^2\lesssim
    T^{1/2}\norm{u}_{X_T}^2.
\end{align*}

Now, we turn to the case $1\ll N_1\ll N_2$. From the
Littlewood-Paley theorem, we have for $p\gs 4$
\begin{align}\label{eqprop13}
    &\norm{\sum_{1\ll N_1\ll N_2}P_{N_1}uP_{N_2}\bar{u}_x}_{L_x^p
    L_T^2}=\norm{\sum_{N_1\gg 1}P_{N_1}u P_{\gg
    N_1}\bar{u}_x}_{L_x^p
    L_T^2}\no\\
    \lesssim&  \sum_{N_1\gg 1}\norm{P_{N_1}u P_{\gg
    N_1}\bar{u}_x}_{L_x^p L_T^2}
    \lesssim \sum_{N_1\gg 1}\norm{P_{N_1}u}_{L_x^p
    L_T^\infty}\norm{P_{\gg N_1}u_x}_{L_x^\infty L_T^2}.
\end{align}
Noticing that
\begin{align*}
    \norm{P_{\gg N_1}u_x}_{L_x^\infty L_T^2}
    \sim & \norm{\left(\sum_M\abs{P_M P_{\gg N_1}u_x}^2\right)^{1/2}}_{L_x^\infty L_T^2}
    \lesssim \left(\sum_{M\gg N_1}\norm{P_Mu_x}_{L_x^\infty
    L_T^2}^2\right)^{1/2}\\
    \lesssim & \norm{u}_{X_T},
\end{align*}
and for $N_1\gg 1$, $\eps=1/p$ and $p\gs 4$
\begin{align}
    \norm{P_{N_1}u}_{L_x^p L_T^\infty}
    =&\norm{P_{N_1}(P_{\lesssim  1}u+P_{\gg 1}u)}_{L_x^p L_T^\infty}
    =\norm{P_{N_1}P_{\gg 1}u}_{L_x^p L_T^\infty}\no\\
    = & N_1^{-\eps} N_1^\eps \norm{P_{N_1}P_{\gg 1}u}_{L_x^p L_T^\infty}
    \sim N_1^{-\eps}  \norm{D_x^\eps P_{N_1}P_{\gg 1}u}_{L_x^p
    L_T^\infty}\no\\
    \lesssim &N_1^{-\eps} \norm{\dx^\frac{1}{4} P_{N_1}P_{\gg
    1}u}_{L_x^4 L_T^\infty}, \label{eqprop140}
\end{align}
we can bound \eqref{eqprop13} by
\begin{align}\label{eqprop14}
    \lesssim &\norm{u}_{X_T} \sum_{N_1\gg 1}N_1^{-\eps}\norm{\dx^\frac{1}{4}
    P_{N_1}P_{\gg
    1}u}_{L_x^4 L_T^\infty}\no\\
    \lesssim & \norm{u}_{X_T} \left(\sum_{N_1\gg 1}N_1^{-2\eps}\right)^{\frac{1}{2}}
    \left(\sum_{N_1\gg 1}\norm{\dx^\frac{1}{4} P_{N_1}P_{\gg
    1}u}_{L_x^4 L_T^\infty}^2\right)^{\frac{1}{2}}\no\\
    \lesssim & \norm{u}_{X_T}\norm{P_{\gg 1}u}_{X_T},
\end{align}
in view of the H\"older inequality. Thus, we have obtained
\begin{align}\label{eqprop141}
    \norm{\sum_{1\ll N_1\ll N_2}P_{N_1}uP_{N_2}\bar{u}_x}_{L_x^p
    L_T^2}\lesssim \norm{u}_{X_T}\norm{P_{\gg 1}u}_{X_T}, \quad
    \forall p\gs 4.
\end{align}
Therefore, we have the desired result \eqref{eqprop1} for any real
number $p\gs 4$.
\end{proof}

\section{Nonlinear estimates}

To state the estimates for the nonlinearities $I_{N,j}$, we define
the function space $Y_T$ equipped with the following norm:
\begin{align*}
    \norm{u}_{Y_T}=\norm{u}_{L_T^\infty H_x^{1/2}}+\norm{\partial_x
    u}_{L_x^\infty L_T^2}+\norm{u}_{L_x^2 L_T^\infty}+\norm{\dx^{\frac{1}{4}}
    u}_{L_x^4 L_T^\infty}.
\end{align*}

We have the following proposition for the nonlinearities.

\begin{proposition}\label{prop2}
Let $u$ be a $H^\infty$-solution to \eqref{dnls}-\eqref{data}. Then,
\begin{align*}
    &\left(\sum_{N\gg 1}\norm{\int_0^t
    S(t-\tau)\sum_{j=1}^5I_{N,j}(\tau)d\tau}_{Y_T}^2\right)^{1/2}\\
    \lesssim &(1+\norm{u}_{X_T}^k)\left[T^{\frac{1}{2}}\norm{u}_{X_T}^{\tilde{k}+1}\norm{P_{\gg
1}u}_{X_T}^{k-\tilde{k}}
    +(1+T^{\frac{1}{4}}\norm{u}_{X_T})\norm{u}_{X_T}^{\tilde{k}}\norm{P_{\gg
    1}u}_{X_T}^{k+1-\tilde{k}}\right]\\
    &+T^{\frac{1}{2}}\left(\norm{u}_{X_T}^{2k-1}
    +\norm{u}_{X_T}^{(5k-2)/2}\right)\norm{P_{\gg
    1}u}_{X_T}\\
    &+T^{\frac{1}{4}}(1+T^{\frac{1}{4}}\norm{u}_{X_T})^{\frac{1}{2}}\norm{u}_{X_T}^{2k-1}\norm{P_{\gg
    1}u}_{X_T}^{\frac{3}{2}}\\
    &+(1+\norm{u}_{X_T}^k)\left[T^{\frac{1}{2}} \norm{u}_{X_T}^k\norm{P_{\gg
    1}u}_{X_T}+(1+T^{\frac{1}{4}}\norm{u}_{X_T})^2\norm{u}_{X_T}^{k-1}\norm{P_{\gg
    1}u}_{X_T}^2\right]\\
    &+T^{\frac{1}{2}}\norm{u}_{X_T}^{3k}\norm{P_{\gg 1}
    u}_{X_T},
\end{align*}
where $\tilde{k}$ denotes the maximal integer that is less than $k$
(i.e. $\tilde{k}=[k]$ if $k$ is not an integer and $\tilde{k}=k-1$
if $k$ is an integer where $[k]$ denotes the maximal integer that is
less than or equal to $k$).
\end{proposition}

We consider each nonlinearity separately.

\subsection{Nonlinear estimates of $I_{N,1}$}

Noting that the term $P_N(\abs{\pn}^k u_x)$ has Fourier support in
$\abs{\xi}\sim N$, we have
\begin{align}
    &P_N(\abs{u}^ku_x)-\abs{\pn}^k\pn[]_x\no\\
    =&P_N((\abs{u}^k-\abs{\pn}^k)u_x)+P_N(\abs{\pn}^k u_x)
    -\abs{\pn}^k P_Nu_x\no\\
    =&P_N((\abs{u}^k-\abs{\pn}^k)u_x)+P_N(\abs{\pn}^k \tilde{P}_N u_x)
    -\abs{\pn}^k P_N\tilde{P}_Nu_x,\label{eqin11}
\end{align}
where $\tilde{P}_N=P_{N/2}+P_N+P_{2N}$.

For the second term in \eqref{eqin11}, we have the following
estimate.

\begin{lemma}\label{lem2}
Let $u$ be a solution of \eqref{dnls}-\eqref{data}. Then, we have
for any $k\gs 4$
\begin{align}\label{eqin18}
    &\left(\sum_{N\gg 1}\norm{P_N(\abs{\pn}^k \tilde{P}_N u_x)
    -\abs{\pn}^k P_N\tilde{P}_Nu_x}_{L_x^1
    L_T^2}^2\right)^{1/2}\\
    \lesssim &T^{\frac{1}{2}} \norm{u}_{X_T}^k \norm{P_{\gg
    1}u}_{X_T}
    +(1+T^{\frac{1}{4}}\norm{u}_{X_T})\norm{u}_{X_T}^{k-1}\norm{P_{\gg
    1}u}_{X_T}^2.\no
\end{align}
\end{lemma}

\begin{proof}
To shift a derivative from the high-frequency function $P_Nu_x$ to
the low-frequency function $\abs{\pn}^k$, we require the following
Leibniz rule for $P_N$ from \cite{KT05}:
\begin{align}\label{eqin12}
(P_N(fg)-fP_N g)(x)=\int_0^1 \left(\int\check{\varphi}_N(y)y
f_x(x-\eta y)g(x-y) dy\right)d\eta.
\end{align}
Thus, we have
\begin{align}\label{eqin14}
    &\norm{P_N(\abs{\pn}^k \tilde{P}_N u_x)-\abs{\pn}^k
    P_N\tilde{P}_Nu_x}_{L_x^1 L_T^2}\no\\
    \lesssim &
    \norm{\check{\varphi}_N(y)y}_{L_y^1}\norm{(\abs{\pn}^k)_x}_{L_{x,T}^2}\norm{\tilde{P}_Nu_x}_{L_x^2
    L_T^\infty}\no\\
    \lesssim & N^{-1}\norm{\check{\varphi}_1(y)y}_{L_y^1}\norm{(\abs{\pn}^k)_x}_{L_{x,T}^2}\norm{\tilde{P}_Nu_x}_{L_x^2
    L_T^\infty}\no\\
    \lesssim &\norm{(\abs{\pn}^k)_x}_{L_{x,T}^2}\norm{\tilde{P}_Nu}_{L_x^2
    L_T^\infty}\no\\
    \lesssim &
    \norm{\pn}_{L_{x}^{2k}L_T^\infty}^{k-2}\norm{\pn_x\overline{\pn}}_{L_x^{k}L_T^2}\norm{\tilde{P}_Nu}_{L_x^2
    L_T^\infty}.
\end{align}
Decomposing $\pn=P_{\ls 1}u+P_{1<\cdot\ll N}u$ for $N\gg 1$, we have
\begin{align}\label{eqin13}
    \pn_x\overline{\pn}=P_{\ls 1}u_x\overline{\pn}+P_{1<\cdot\ll
    N}u_x\overline{P_{\ls 1}u}+P_{1<\cdot\ll N}u_x\overline{P_{1<\cdot\ll N}u}.
\end{align}

For the first term in \eqref{eqin13}, we have
\begin{align*}
    \norm{P_{\ls 1}u_x\overline{\pn}}_{L_x^{k}L_T^2}\lesssim
    &\norm{P_{\ls
    1}u_x}_{L_x^{2k}L_T^2}\norm{\pn}_{L_{x}^{2k}L_T^\infty}\\
    \lesssim &T^{1/2}\norm{P_{\ls
    1}u}_{L_x^{2k}L_T^\infty}\norm{\pn}_{L_{x}^{2k}L_T^\infty}.
\end{align*}
By the Littlewood-Paley theorem, we can obtain
\begin{align}\label{eq5.50}
    \norm{\pn}_{L_{x}^{2k}L_T^\infty}\sim &\norm{\left(\sum_M \abs{P_M
    \pn}^2\right)^{1/2}}_{L_{x}^{2k}L_T^\infty}
    \lesssim \left(\sum_M \norm{P_M
    \pn}_{L_{x}^{2k}L_T^\infty}^2\right)^{1/2}\no\\
    \lesssim &\left(\sum_{M\ll N} \norm{P_M
    u}_{L_{x}^{2k}L_T^\infty}^2\right)^{1/2}
    \lesssim \left(\sum_{M\gg 1} \norm{P_M
    u}_{L_{x}^{2k}L_T^\infty}^2\right)^{1/2}\no\\
    \lesssim &\left(\sum_{M\gg 1} \norm{\langle D_x\rangle^{\frac{1}{4}} P_M
    u}_{L_{x}^4 L_T^\infty}^2\right)^{1/2}\lesssim \norm{P_{\gg  1}u}_{X_T}.
\end{align}
In the similar way, we have
\begin{align*}
    \norm{P_{\ls 1}u}_{L_x^{2k}L_T^\infty}\lesssim \norm{u}_{X_T}.
\end{align*}
Thus,
\begin{align}\label{eqin15}
    \norm{P_{\ls 1}u_x\overline{\pn}}_{L_x^{k}L_T^2}\lesssim T^{1/2}\norm{u}_{X_T}\norm{P_{\gg  1}u}_{X_T}.
\end{align}

For the last two term in \eqref{eqin13}, in a similar way as in the
proof of Proposition \ref{prop1}, we can obtain the following bound:
\begin{align}
    \norm{P_{\ls 1}\bar{u} P_{1\ll \cdot\ll N}u_x}_{L_x^k
    L_T^2}\lesssim &T^{\frac{1}{2}} \norm{u}_{X_T}^2+\norm{u}_{X_T}\norm{P_{\gg
    1}u}_{X_T},\label{eqin16}\\
    \norm{P_{1\ll \cdot\ll N}\bar{u} P_{1\ll \cdot\ll N}u_x}_{L_x^k
    L_T^2}\lesssim &T^{\frac{1}{2}} \norm{u}_{X_T}^2+(1+T^{\frac{1}{4}}\norm{u}_{X_T})\norm{u}_{X_T}\norm{P_{\gg
    1}u}_{X_T}.\label{eqin17}
\end{align}
From the Sobolev embedding theorem and
\eqref{eqin15}-\eqref{eqin17}, we obtain that \eqref{eqin14} can be
bounded by
\begin{align*}
    \lesssim \left(T^{\frac{1}{2}} \norm{u}_{X_T}^k+(1+T^{\frac{1}{4}}\norm{u}_{X_T})\norm{u}_{X_T}^{k-1}\norm{P_{\gg
    1}u}_{X_T}\right)\norm{\tilde{P}_Nu}_{L_x^2
    L_T^\infty}.
\end{align*}
Thus, we can bound \eqref{eqin18} by
\begin{align*}
    \lesssim &\left(T^{\frac{1}{2}} \norm{u}_{X_T}^k+(1+T^{\frac{1}{4}}\norm{u}_{X_T})\norm{u}_{X_T}^{k-1}\norm{P_{\gg
    1}u}_{X_T}\right)\left(\sum_{N\gg 1}\norm{P_N u}_{L_x^2
    L_T^\infty}^2\right)^{1/2}\\
    \lesssim &T^{\frac{1}{2}} \norm{u}_{X_T}^k \norm{P_{\gg
    1}u}_{X_T}
    +(1+T^{\frac{1}{4}}\norm{u}_{X_T})\norm{u}_{X_T}^{k-1}\norm{P_{\gg
    1}u}_{X_T}^2,
\end{align*}
which yields the desired result.
\end{proof}

For the first term in \eqref{eqin11}, we have the following
estimate:

\begin{lemma}\label{lem3}
Let $u$ be a solution of \eqref{dnls}-\eqref{data}. Then, we have
for any $k\gs 4$
\begin{align}\label{eqin206}
&\left(\sum_{N\gg 1}\norm{P_N((|u|^k-|P_{\ll N}u|^k)u_x)}_{L_x^1
L_T^2}^2\right)^{1/2}\no\\
\lesssim &T^{\frac{1}{2}}\norm{u}_{X_T}^{\tilde{k}+1}\norm{P_{\gg
1}u}_{X_T}^{k-\tilde{k}}
    +(1+T^{\frac{1}{4}}\norm{u}_{X_T})\norm{u}_{X_T}^{\tilde{k}}\norm{P_{\gg
    1}u}_{X_T}^{k+1-\tilde{k}}.
\end{align}
\end{lemma}

\begin{proof}
We split \eqref{eqin11} into several terms for $N\gg 1$ and $k\gs 4$
\begin{align}
    &P_N((|u|^k-|P_{\ll N}u|^k)u_x)\label{eqin201}\\
    =&P_N(|u|^{k-2}\bar{u}u_x P_{\gtrsim N}u )\label{eqin202}\\
    &+P_N((|u|^{k-2}-|P_{\ll N}u|^{k-2})\bar{u}u_xP_{\ll N}u)\label{eqin203}\\
    &+P_N(|P_{\ll
    N}u|^{k-2}u_xP_{\ll N}uP_{\gtrsim N}\bar{u}).\label{eqin204}
\end{align}
Notice that
\begin{align*}
    \norm{P_{\gtrsim N}u}_{L_x^k L_T^\infty}
    \lesssim& \left(\sum_{M\gtrsim
    N}\norm{P_Mu}_{L_x^k L_T^\infty}^2\right)^{1/2}
    \lesssim \left(\sum_{M\gtrsim N}N^{-2\eps_k}M^{2\eps_k}\norm{P_Mu}_{L_x^k
    L_T^\infty}^2\right)^{1/2}\\
    \lesssim &\left(\sum_{M\gtrsim N}N^{-2\eps_k}\norm{D_x^{\eps_k} P_Mu}_{L_x^k
    L_T^\infty}^2\right)^{1/2}\\
    \lesssim &N^{-\eps_k}\left(\sum_{M\gtrsim N}\norm{\dx^{\frac{1}{4}}
    P_Mu}_{L_x^4 L_T^\infty}^2\right)^{1/2}\\
    \lesssim &N^{-\eps_k}\norm{P_{\gg  1}u}_{X_T}, \quad \forall k\gs
    4,
\end{align*}
where $\eps_k>0$ is defined by $\eps_k=1/k$.

Thus, for the first term \eqref{eqin202}, from the fact
$\norm{\check{\varphi}_N}_{L^1}\lesssim 1$ and
Proposition~\ref{prop1}, we have for $k\gs 4$
\begin{align*}
    &\norm{P_N(|u|^{k-2}\bar{u}u_x P_{\gtrsim N}u )}_{L_x^1 L_T^2}
    \lesssim \norm{|u|^{k-2}\bar{u}u_x P_{\gtrsim N}u }_{L_x^1
    L_T^2}\\
    \lesssim & \norm{u}_{L_x^{k}L_T^\infty}^{k-2}\norm{\bar{u}u_x}_{L_x^k L_T^2}\norm{P_{\gtrsim
    N}u}_{L_x^k L_T^\infty}\\
    \lesssim & N^{-\eps_k}\left[T^{\frac{1}{2}}\norm{u}_{X_T}^k\norm{P_{\gg
    1}u}_{X_T}
    +(1+T^{\frac{1}{4}}\norm{u}_{X_T})\norm{u}_{X_T}^k\norm{P_{\gg
    1}u}_{X_T}^2\right].
\end{align*}
Therefore, we obtain, for any $k\gs 4$, that
\begin{align*}
    &\left(\sum_{N\gg 1}\norm{P_N(|u|^{k-2}\bar{u}u_x P_{\gtrsim N}u )}_{L_x^1
    L_T^2}^2\right)^{1/2}\\
    \lesssim & T^{\frac{1}{2}}\norm{u}_{X_T}^k\norm{P_{\gg
    1}u}_{X_T} +(1+T^{\frac{1}{4}}\norm{u}_{X_T})\norm{u}_{X_T}^k\norm{P_{\gg
    1}u}_{X_T}^2.
\end{align*}

For \eqref{eqin204}, in the same way as the case \eqref{eqin202}, we
have
\begin{align*}
    &\left(\sum_{N\gg 1}\norm{P_N(|P_{\ll
    N}u|^{k-2}u_xP_{\ll N}uP_{\gtrsim N}\bar{u})}_{L_x^1
    L_T^2}^2\right)^{1/2}\\
    \lesssim & T^{\frac{1}{2}}\norm{u}_{X_T}^k\norm{P_{\gg
    1}u}_{X_T} +(1+T^{\frac{1}{4}}\norm{u}_{X_T})\norm{u}_{X_T}^k\norm{P_{\gg
    1}u}_{X_T}^2.
\end{align*}

Now, we derive the estimate for \eqref{eqin203} by using the
induction argument in $k$.

For $k=4$, we have
\begin{align*}
    \abs{|u|^{k-2}-|P_{\ll N}u|^{k-2}}\lesssim \abs{P_{\gtrsim N}u}^{2}+\abs{P_{\gtrsim N}u \overline{P_{\ll N}u}}.
\end{align*}
From the Young inequality, the H\"older inequality, \eqref{eq5.50}
and Proposition~\ref{prop1}, we can get for $k=4$
\begin{align*}
    &\norm{ P_N((|u|^{k-2}-|P_{\ll N}u|^{k-2})\bar{u}u_xP_{\ll
    N}u)}_{L_x^1 L_T^2}\\
    \lesssim & \norm{|u|^{k-2}-|P_{\ll
    N}u|^{k-2}}_{L_x^{2k/(k-2)}L_T^\infty}\norm{\bar{u}u_x}_{L_x^k L_T^2}\norm{P_{\ll
    N}u}_{L_x^2 L_T^\infty}\\
    \lesssim & (\norm{P_{\gtrsim
    N}u}_{L_x^{8}L_T^\infty}^{2}+\norm{P_{\gtrsim N}u }_{L_x^{8}L_T^\infty}\norm{P_{\ll N}u}_{L_x^{8}L_T^\infty})
    \norm{\bar{u}u_x}_{L_x^k L_T^2}\norm{P_{\ll
    N}u}_{L_x^2 L_T^\infty}\\
    \lesssim & N^{-1/8}\norm{P_{\gg  1}u}_{X_T}^{2}\norm{\bar{u}u_x}_{L_x^k L_T^2}\norm{P_{\ll
    N}u}_{L_x^2 L_T^\infty}\\
    \lesssim & N^{-1/8}\left[T^{\frac{1}{2}}\norm{u}_{X_T}^3\norm{P_{\gg  1}u}_{X_T}^{k-2}
    +(1+T^{\frac{1}{4}}\norm{u}_{X_T})\norm{u}_{X_T}^2\norm{P_{\gg
    1}u}_{X_T}^{k-1}\right].
\end{align*}

From the triangle inequality for complex number, i.e.
$\abs{|z_1|-|z_2|}\ls |z_1-z_2|$ for $z_1,z_2\in \mathbb{C}$, we can
get $\abs{|z_1|^\theta-|z_2|^\theta}\ls |z_1-z_2|^\theta$ for any
$\theta\in (0,1]$.

For $k\in(4,5]$, we have
\begin{align*}
    &\abs{|u|^{k-2}-|P_{\ll N}u|^{k-2}}\\
    \lesssim
    &|u|^{2}\abs{|u|^{k-4}-|P_{\ll N}u|^{k-4}}+|P_{\ll N}u|^{k-4}\abs{\abs{u}^2-|P_{\ll
    N}u|^2}\\
    \lesssim &|u|^{2}\abs{P_{\gtrsim N}u}^{k-4}+|P_{\ll N}u|^{k-4}
    \abs{P_{\gtrsim N}u}^{2}+\abs{P_{\gtrsim N}u} \abs{P_{\ll
    N}u}^{k-3}.
\end{align*}
Then
\begin{align*}
    &\norm{ P_N((|u|^{k-2}-|P_{\ll N}u|^{k-2})\bar{u}u_xP_{\ll
    N}u)}_{L_x^1 L_T^2}\\
    \lesssim & \left[\norm{u}_{L_x^{2k}L_T^\infty}^{2}\norm{P_{\gtrsim
    N}u}_{L_x^{2k}L_T^\infty}^{k-4}+\norm{P_{\ll N}u}_{L_x^{2k}L_T^\infty}^{k-4}
    \norm{P_{\gtrsim N}u}_{L_x^{2k}L_T^\infty}^{2}\right.\\
    &\qquad \left.+\norm{P_{\gtrsim N}u}_{L_x^{2k}L_T^\infty}\norm{P_{\ll
    N}u}_{L_x^{2k}L_T^\infty}^{k-3}\right]\norm{\bar{u}u_x}_{L_x^k L_T^2}\norm{P_{\ll
    N}u}_{L_x^{2}L_T^\infty}\\
    \lesssim &N^{-\eps_k}\left[T^{\frac{1}{2}}\norm{u}_{X_T}^4\norm{P_{\gg  1}u}_{X_T}^{k-3}
    +(1+T^{\frac{1}{4}}\norm{u}_{X_T})\norm{u}_{X_T}^3\norm{P_{\gg
    1}u}_{X_T}^{k-2}\right],
\end{align*}
where $\eps_k=(k-4)/2k$ for $k\in(4,5]$. By the same procedure, we
can obtain for any $k\gs 4$
\begin{align*}
    &\norm{ P_N((|u|^{k-2}-|P_{\ll N}u|^{k-2})\bar{u}u_xP_{\ll
    N}u)}_{L_x^1 L_T^2}\\
    \lesssim &N^{-\eps_k}\left[T^{\frac{1}{2}}\norm{u}_{X_T}^{\tilde{k}+1}\norm{P_{\gg  1}u}_{X_T}^{k-\tilde{k}}
    +(1+T^{\frac{1}{4}}\norm{u}_{X_T})\norm{u}_{X_T}^{\tilde{k}}\norm{P_{\gg
    1}u}_{X_T}^{k+1-\tilde{k}}\right],
\end{align*}
where $\eps_k=(k-\tilde{k})/2k>0$.  Therefore, we have for any $k\gs
4$
\begin{align}\label{eqin205}
    &\left(\sum_{N\gg 1}\norm{ P_N((|u|^{k-2}-|P_{\ll N}u|^{k-2})\bar{u}u_xP_{\ll
    N}u)}_{L_x^1 L_T^2}^2\right)^{1/2}\no\\
    \lesssim &T^{\frac{1}{2}}\norm{u}_{X_T}^{\tilde{k}+1}\norm{P_{\gg  1}u}_{X_T}^{k-\tilde{k}}
    +(1+T^{\frac{1}{4}}\norm{u}_{X_T})\norm{u}_{X_T}^{\tilde{k}}\norm{P_{\gg
    1}u}_{X_T}^{k+1-\tilde{k}}.
\end{align}
Thus, we have proved the desired result.
\end{proof}

\begin{remark}
From the proof of Lemma~\ref{lem3}, we can see that
\begin{align}\label{eqin207}
&\left(\sum_{N\gg 1}\norm{P_N((|u|^k-|P_{\ll
N}u|^k)u_x)}_{L_x^{\frac{1}{1-\eps}}
L_T^2}^2\right)^{1/2}\no\\
\lesssim &T^{\frac{1}{2}}\norm{u}_{X_T}^{\tilde{k}+1}\norm{P_{\gg
1}u}_{X_T}^{k-\tilde{k}}
    +(1+T^{\frac{1}{4}}\norm{u}_{X_T})\norm{u}_{X_T}^{\tilde{k}}\norm{P_{\gg
    1}u}_{X_T}^{k+1-\tilde{k}},
\end{align}
holds for any $\eps\in[0,1)$ in view of Proposition~\ref{prop1}.
\end{remark}

We turn to the proof of Proposition~\ref{prop2} for the nonlinearity
$I_{N,1}$. We also consider the decomposition in \eqref{eqin11}. For
convenience, we denote $B_N=P_N(|P_{\ll
N}u|^k\tilde{P}_Nu_x)-|P_{\ll N}u|^kP_N\tilde{P}_Nu_x$. From
\eqref{eq5}, \eqref{eq12}, \eqref{eq7} and \eqref{eq8},
 we have
\begin{align}
    &\left(\sum_{N\gg 1}\norm{\int_0^t S(t-\tau)\ep B_N
    d\tau}_{Y_T}^2\right)^{1/2}\no\\
    \lesssim &\left(\sum_{N\gg 1}\norm{ B_N
    }_{L_x^1 L_T^2}^2\right)^{1/2}+\left(\sum_{N\gg 1}\left(\sum_{M}\norm{P_M(\ep
    B_N)}_{L_x^1 L_T^2}\right)^2\right)^{1/2}.\label{eq15}
\end{align}
By Lemma~\ref{lem2}, the first term can be bounded by
\begin{align*}
\lesssim T^{\frac{1}{2}} \norm{u}_{X_T}^k \norm{P_{\gg
    1}u}_{X_T}
    +(1+T^{\frac{1}{4}}\norm{u}_{X_T})\norm{u}_{X_T}^{k-1}\norm{P_{\gg
    1}u}_{X_T}^2.
\end{align*}
For the second term, we split the sum $\sum_M$ into three parts
$\sum_{M\sim N}+\sum_{M\ll N}$ $+\sum_{M\gg N}$ as in \cite{KT05}.
For the part of $M\sim N$, it is the same as Lemma~\ref{lem2} by
summing in $M$ such that $M\sim N$. For the part $M\ll N$, we can
add the projection operator $P_{\sim N}$ to $\ep$ since $B_N$ has
Fourier support in $\abs{\xi}\sim N$. Thus, by the H\"older
inequality, we have
\begin{align}
    &\left(\sum_{N\gg 1}\left(\sum_{M\ll N}\norm{P_M(\ep
    B_N)}_{L_x^1 L_T^2}\right)^2\right)^{1/2}\no\\
    \lesssim &\left(\sum_{N\gg 1}\left(\sum_{M\ll N}\norm{P_{\sim N}\ep
    B_N}_{L_x^1 L_T^2}\right)^2\right)^{1/2}\no\\
    \lesssim &\left(\sum_{N\gg 1}(\ln N)^2\norm{P_{\sim
    N}\ep}_{L_{x}^{1/\eps}L_T^\infty}^2
    \norm{B_N}_{L_x^{1/(1-\eps)}
    L_T^{2}}^2\right)^{1/2},\label{eq13}
\end{align}
where $\eps\in(0,1/k)$.

 By the Bernstein inequality, we have
\begin{align*}
    &N\norm{P_{\sim N}\ep}_{L_{x}^{1/\eps}L_T^\infty}\lesssim \norm{\partial_x P_{\sim
    N}\ep}_{L_{x}^{1/\eps}L_T^\infty}\\
    \lesssim &\norm{P_{\ll N}u}_{L_{x}^{k/\eps}L_T^\infty}^k\lesssim
    \norm{u}_{X_T}^k,
\end{align*}
and from \eqref{eqin12} and the H\"older inequality, we can get, as
a similar way as in \eqref{eqin14}, that
\begin{align*}
    &\norm{B_N}_{L_x^{1/(1-\eps)} L_T^{2}}
    =\norm{(|P_{\ll N}u|^k)_x}_{L_x^{2/(1-2\eps)} L_T^2}\norm{\tilde{P}_N u}_{L_x^2
    L_T^\infty}\\
    \lesssim &\norm{P_{\ll N}u}_{L_x^{\frac{2k(k-2)}{k-2-2k\eps}} L_T^\infty}^{k-2}\norm{P_{\ll N}u_xP_{\ll
    N}\bar{u}}_{L_x^{k}L_T^2}\norm{\tilde{P}_N u}_{L_x^2
    L_T^\infty}\\
    \lesssim &\left(T^{\frac{1}{2}} \norm{u}_{X_T}^k+(1+T^{\frac{1}{4}}\norm{u}_{X_T})\norm{u}_{X_T}^{k-1}\norm{P_{\gg
    1}u}_{X_T}\right)\norm{P_N u}_{L_x^2 L_T^\infty}.
\end{align*}
Thus, \eqref{eq13} can be bounded by
\begin{align*}
   \lesssim & \norm{u}_{X_T}^k\left(T^{\frac{1}{2}} \norm{u}_{X_T}^k+(1+T^{\frac{1}{4}}
   \norm{u}_{X_T})\norm{u}_{X_T}^{k-1}\norm{P_{\gg
    1}u}_{X_T}\right) \left(\sum_{N\gg 1} \norm{P_N u}_{L_x^2
    L_T^\infty}^2\right)^{1/2}\\
    \lesssim &T^{\frac{1}{2}} \norm{u}_{X_T}^{2k}\norm{P_{\gg
    1}u}_{X_T}+(1+T^{\frac{1}{4}}
   \norm{u}_{X_T})\norm{u}_{X_T}^{2k-1}\norm{P_{\gg
    1}u}_{X_T}^2.
\end{align*}
For the part $M\gg N$, we can add the projection operator $P_M$ to
$\ep$. In a similar way with the part $M\ll N$, we have
\begin{align*}
    &\left(\sum_{N\gg 1}\left(\sum_{M\gg N}\norm{P_M(\ep
    B_N)}_{L_x^1 L_T^2}\right)^2\right)^{1/2}\\
    \lesssim &\left(\sum_{N\gg 1}\left(\sum_{M\gg N}\norm{P_M \ep}_{L_{x}^{1/\eps}L_T^\infty}
    \norm{B_N}_{L_x^{1/(1-\eps)}
    L_T^{2}}\right)^2\right)^{1/2}\\
    \lesssim &\Bigg(\sum_{N\gg 1}\Big(\sum_{M\gg N}\frac{1}{M}\Big(T^{\frac{1}{2}}
    \norm{u}_{X_T}^{2k}\norm{P_N u}_{L_x^2 L_T^\infty}\\
    &\qquad\qquad\qquad+(1+T^{\frac{1}{4}}\norm{u}_{X_T})\norm{u}_{X_T}^{2k-1}\norm{P_{\gg
    1}u}_{X_T}\norm{P_N u}_{L_x^2 L_T^\infty}\Big)\Big)^2\Bigg)^{1/2}\\
    \lesssim &T^{\frac{1}{2}} \norm{u}_{X_T}^{2k}\norm{P_{\gg
    1}u}_{X_T}+(1+T^{\frac{1}{4}}
   \norm{u}_{X_T})\norm{u}_{X_T}^{2k-1}\norm{P_{\gg
    1}u}_{X_T}^2.
\end{align*}
For the first term in \eqref{eqin11}, we denote it by $A_N$, i.e.
$A_N=P_N((|u|^k-|P_{\ll N}u|^k)u_x)$. Similarly, from \eqref{eq5},
\eqref{eq12}, \eqref{eq7} and \eqref{eq8},
 we can get
\begin{align}
    &\left(\sum_{N\gg 1}\norm{\int_0^t S(t-\tau)\ep A_N
    d\tau}_{Y_T}^2\right)^{1/2}\no\\
    \lesssim &\left(\sum_{N\gg 1}\norm{ A_N
    }_{L_x^1 L_T^2}^2\right)^{1/2}+\left(\sum_{N\gg 1}\left(\sum_{M}\norm{P_M(\ep
    A_N)}_{L_x^1 L_T^2}\right)^2\right)^{1/2}.\label{eq14}
\end{align}
From Lemma~\ref{lem3}, the first term is bounded by
\begin{align*}
    \lesssim T^{\frac{1}{2}}\norm{u}_{X_T}^{\tilde{k}+1}\norm{P_{\gg
1}u}_{X_T}^{k-\tilde{k}}
    +(1+T^{\frac{1}{4}}\norm{u}_{X_T})\norm{u}_{X_T}^{\tilde{k}}\norm{P_{\gg
    1}u}_{X_T}^{k+1-\tilde{k}}.
\end{align*}
Noticing that \eqref{eqin207}, and in the same way as in dealing
with the second term of \eqref{eq15}, we can bound the second term
of \eqref{eq14} by
\begin{align*}
    \lesssim T^{\frac{1}{2}}\norm{u}_{X_T}^{k+\tilde{k}+1}\norm{P_{\gg
1}u}_{X_T}^{k-\tilde{k}}
    +(1+T^{\frac{1}{4}}\norm{u}_{X_T})\norm{u}_{X_T}^{k+\tilde{k}}\norm{P_{\gg
    1}u}_{X_T}^{k+1-\tilde{k}}.
\end{align*}
Therefore, we have obtained
\begin{align*}
    &\left(\sum_{N\gg 1}\norm{\int_0^t
    S(t-\tau)I_{N,1}(\tau)d\tau}_{Y_T}^2\right)^{1/2}\\
    \lesssim &(1+\norm{u}_{X_T}^k)\left[T^{\frac{1}{2}}\norm{u}_{X_T}^{\tilde{k}+1}\norm{P_{\gg
1}u}_{X_T}^{k-\tilde{k}}
    +(1+T^{\frac{1}{4}}\norm{u}_{X_T})\norm{u}_{X_T}^{\tilde{k}}\norm{P_{\gg
    1}u}_{X_T}^{k+1-\tilde{k}}\right].
\end{align*}

\subsection{Nonlinear estimates of $I_{N,2}$}

From \eqref{eq4}, \eqref{eq9}, \eqref{eq10} and \eqref{eq11}, we
have
\begin{align}
    &\left(\sum_{N\gg 1}\norm{\int_0^t
    S(t-\tau)I_{N,2}(\tau)d\tau}_{Y_T}^2\right)^{1/2}\no\\
    \lesssim &\left(\sum_{N\gg 1}\norm{\ep P_N u
    B_{N,2}}_{L_T^1H_x^{1/2}}^2\right)^{1/2}\label{eq16}\\
    &+\left(\sum_{N\gg 1}\left(\sum_{M}\norm{P_M(\ep P_N u
    B_{N,2})}_{L_T^1H_x^{1/2}}\right)^2\right)^{1/2},\label{eq17}
\end{align}
where $B_{N,2}=\int_{-\infty}^x
\abs{\pn}^{k-4}\left[(\overline{\pn_x}\pn)^2-(\pn_x\overline{\pn})^2\right]dy$.

For the first term \eqref{eq16}, from Lemma~\ref{lem1} and the
H\"older inequality, it can be bounded by
\begin{align}
    \lesssim &\left(\sum_{N\gg 1}\norm{P_N u
    B_{N,2}}_{L_T^1L_x^2}^2\right.\no\\
    &\left.\quad +\norm{P_N u
    B_{N,2}}_{L_T^1L_x^2}\norm{\partial_x(\ep P_N u
    B_{N,2})}_{L_T^1L_x^2}\right)^{1/2}\no\\
    \lesssim &\Big(\sum_{N\gg 1}\norm{P_N u}_{L_T^\infty L_x^2}^2
    \norm{B_{N,2}}_{L_T^1L_x^\infty}^2+\norm{P_N u}_{L_T^\infty H_x^{1/2}}^2
    \norm{B_{N,2}}_{L_T^1L_x^\infty}^2
    \no\\
    &+\norm{P_N u}_{L_T^\infty L_x^2}
    \norm{B_{N,2}}_{L_T^1L_x^\infty}^2\norm{\pn}_{L_x^{2(k+1)}L_T^\infty}^k\norm{P_Nu}_{L_x^{2(k+1)}L_T^\infty}\no\\
    &+\norm{P_N u}_{L_T^\infty L_x^2}
    \norm{B_{N,2}}_{L_T^1L_x^\infty}\norm{P_N u}_{L_x^\infty L_T^2}
    \norm{\partial_x B_{N,2}}_{L_{x,T}^2}\Big)^{1/2}.\label{eq18}
\end{align}
By the H\"older inequality, we have for $k\gs 5$
\begin{align}
  \norm{B_{N,2}}_{L_T^1L_x^\infty}
    \lesssim   & \norm{\abs{\pn}^{k-4}
    \left[(\overline{\pn_x}\pn)^2-(\pn_x\overline{\pn})^2\right]}_{L_{x,T}^1}\no\\
    \lesssim
    &\norm{\overline{\pn}\pn_x}_{L_x^4L_T^2}^2\norm{\pn}_{L_x^{2(k-4)}L_T^\infty}^{k-4}\no\\
    \lesssim &T \norm{u}_{X_T}^{k}+(1+T^{\frac{1}{4}}\norm{u}_{X_T})^2\norm{u}_{X_T}^{k-2}\norm{P_{\gg
    1}u}_{X_T}^2,\label{eq21}
\end{align}
and from Proposition~\ref{prop1} and the proof of Lemma~\ref{lem2},
\begin{align*}
    &\norm{P_N u}_{L_x^\infty L_T^2}\norm{\partial_x
    B_{N,2}}_{L_{x,T}^2}\\
    = & \norm{P_N u}_{L_x^\infty L_T^2}\norm{\abs{\pn}^{k-4}
    \left[(\overline{\pn_x}\pn)^2-(\pn_x\overline{\pn})^2\right]}_{L_{x,T}^2}\\
    \lesssim & \norm{P_N u_x}_{L_x^\infty L_T^2}\norm{\pn_x P_{\ll
    N}\bar{u}}_{L_x^{2(k-1)}L_T^2}\norm{\pn}_{L_x^{2(k-1)}L_T^\infty}^{k-2}\\
    \lesssim &\norm{P_N u_x}_{L_x^\infty L_T^2}\left[T^{\frac{1}{2}}
    \norm{u}_{X_T}^k+(1+T^{\frac{1}{4}}\norm{u}_{X_T})\norm{u}_{X_T}^{k-1}\norm{P_{\gg
    1}u}_{X_T}\right].
\end{align*}
Thus, we can bound \eqref{eq16} by
\begin{align*}
    \lesssim (1+\norm{u}_{X_T}^k)\left(T^{1/2}\norm{u}_{X_T}^k \norm{P_{\gg
    1}u}_{X_T}+(1+T^{\frac{1}{4}}\norm{u}_{X_T})^2\norm{u}_{X_T}^{k-2}\norm{P_{\gg
    1}u}_{X_T}^3\right).
\end{align*}

For \eqref{eq17}, we split the sum $\sum_M$ into two parts
$\sum_{M\lesssim N}+\sum_{M\gg N}$, which gives the bound by
\begin{align}\label{eq19}
    \lesssim &\left(\sum_{N\gg 1}\left(\sum_{M\lesssim N}\langle
    M\rangle^{\frac{1}{2}}\norm{P_Nu
    B_{N,2}}_{L_T^1L_x^2}\right)^2\right)^{1/2}\\
    &+\left(\sum_{N\gg 1}\left(\sum_{M\gg N}\norm{P_M \dx^{\frac{1}{2}}\ep P_Nu
    B_{N,2}}_{L_T^1L_x^2}\right)^2\right)^{1/2}. \label{eq20}
\end{align}
For the first term \eqref{eq19}, noticing that $\sum_{M\lesssim
N}\langle M\rangle^{1/2}\lesssim N^{1/2}$ and \eqref{eq21}, we can
bound it by
\begin{align*}
    \lesssim & \left(\sum_{N\gg 1}\left(\norm{D_x^{\frac{1}{2}}P_Nu}_{L_T^\infty L_x^2}
    \norm{B_{N,2}}_{L_T^1L_x^\infty}\right)^2\right)^{1/2}\\
    \lesssim &\norm{u}_{X_T}^k\norm{P_{\gg
    1}u}_{X_T}.
\end{align*}
For the second term \eqref{eq20}, in a similar way with
\eqref{eq18}, we bound it by
\begin{align*}
    \lesssim &\left(\sum_{N\gg 1}\left(\sum_{M\gg N}M^{-\frac{1}{2}}\norm{P_M \dx\ep P_Nu
    B_{N,2}}_{L_T^1L_x^2}\right)^2\right)^{1/2}\\
    \lesssim &\left(\sum_{N\gg 1}\norm{P_Nu
    B_{N,2}}_{L_T^1L_x^2}^2+\norm{\partial_x\ep P_Nu
    B_{N,2}}_{L_T^1L_x^2}^2\right)^{1/2}\\
    \lesssim & (1+\norm{u}_{X_T}^k)\left(T^{1/2}\norm{u}_{X_T}^k \norm{P_{\gg
    1}u}_{X_T}+(1+T^{\frac{1}{4}}\norm{u}_{X_T})^2\norm{u}_{X_T}^{k-2}\norm{P_{\gg
    1}u}_{X_T}^3\right).
\end{align*}
Therefore, we obtain
\begin{align*}
    &\left(\sum_{N\gg 1}\norm{\int_0^t
    S(t-\tau)I_{N,2}(\tau)d\tau}_{Y_T}^2\right)^{1/2}\no\\
    \lesssim & (1+\norm{u}_{X_T}^k)\left(T^{1/2}\norm{u}_{X_T}^k \norm{P_{\gg
    1}u}_{X_T}+(1+T^{\frac{1}{4}}\norm{u}_{X_T})^2\norm{u}_{X_T}^{k-2}\norm{P_{\gg
    1}u}_{X_T}^3\right).
\end{align*}

\subsection{Nonlinear estimates of $I_{N,3}$}

From \eqref{eq4}, \eqref{eq9}, \eqref{eq10} and \eqref{eq11}, we
have
\begin{align}
    &\left(\sum_{N\gg 1}\norm{\int_0^t
    S(t-\tau)I_{N,2}(\tau)d\tau}_{Y_T}^2\right)^{1/2}\label{eq24}\\
    \lesssim &\left(\sum_{N\gg 1}\norm{\ep P_N u
    B_{N,3}}_{L_T^1H_x^{1/2}}^2\right)^{1/2}\label{eq22}\\
    &+\left(\sum_{N\gg 1}\left(\sum_{M}\norm{P_M(\ep P_N u
    B_{N,3})}_{L_T^1H_x^{1/2}}\right)^2\right)^{1/2},\label{eq23}
\end{align}
where $B_{N,3}=\int_{-\infty}^x \abs{\pn}^{k-2}P_{\ll
N}|u|^k(u_x+\bar{u}_x)dy$.

By H\"older inequality, we get
\begin{align*}
    \norm{B_{N,3}}_{L_T^1L_x^\infty}\lesssim &\norm{\abs{\pn}^{k-2}P_{\ll
    N}|u|^k(u_x+\bar{u}_x)}_{L_{x,T}^1}\\
    \lesssim
    &T^{\frac{1}{2}}\norm{\pn}_{L_x^{2k-2}L_T^\infty}^{k-2}\norm{P_{\ll
    N}|u|^k(u_x+\bar{u}_x)}_{L_x^{(2k-2)/k}L_T^2}\\
    \lesssim &T^{\frac{1}{2}}\norm{\pn}_{L_x^{2k-2}L_T^\infty}^{k-2}
    \norm{u}_{L_x^{2k-2}L_T^\infty}^k \norm{u_x}_{L_x^\infty
    L_T^2}\\
    \lesssim &T^{\frac{1}{2}}\norm{u}_{X_T}^{2k-1}.
\end{align*}
From the H\"older inequality and Proposition~\ref{prop1}, we have
\begin{align*}
    \norm{\partial_x B_{N,3}}_{L_{x,T}^2}=&\norm{\abs{\pn}^{k-2}P_{\ll
    N}|u|^k(u_x+\bar{u}_x)}_{L_{x,T}^2}\\
    \lesssim
    &\norm{\pn}_{L_x^{4k-4}L_T^\infty}^{k-2}
    \norm{u}_{L_x^{4k-4}L_T^\infty}^k\norm{\bar{u}u_x}_{L_x^{4k-4}L_T^2}\\
    \lesssim &T^{\frac{1}{2}}\norm{u}_{X_T}^{2k}
    +(1+T^{\frac{1}{4}}\norm{u}_{X_T})\norm{u}_{X_T}^{2k-1}\norm{P_{\gg
    1}u}_{X_T}.
\end{align*}
In addition, for $N\gg 1$, we have $\norm{P_Nu}_{L_x^\infty
L_T^2}\lesssim \norm{P_Nu_x}_{L_x^\infty L_T^2}$. Thus, in the same
way as in the case $I_{N,2}$, we can bound \eqref{eq24} by
\begin{align*}
    \lesssim &T^{\frac{1}{2}}\left(\norm{u}_{X_T}^{2k-1}
    +\norm{u}_{X_T}^{(5k-2)/2}\right)\norm{P_{\gg
    1}u}_{X_T}\\
    &+T^{\frac{1}{4}}(1+T^{\frac{1}{4}}\norm{u}_{X_T})^{\frac{1}{2}}\norm{u}_{X_T}^{2k-1}\norm{P_{\gg
    1}u}_{X_T}^{\frac{3}{2}}.
\end{align*}

\subsection{Nonlinear estimates of $I_{N,4}$}

From \eqref{eq5}, \eqref{eq12}, \eqref{eq7} and \eqref{eq8}, we have
\begin{align}
    &\left(\sum_{N\gg 1}\norm{\int_0^t S(t-\tau)\ep B_{N,4}
    d\tau}_{Y_T}^2\right)^{1/2}\no\\
    \lesssim &\left(\sum_{N\gg 1}\norm{ B_{N,4}
    }_{L_x^1 L_T^2}^2\right)^{1/2}\no\\
    &+\left(\sum_{N\gg 1}\left(\sum_{M}\norm{P_M(\ep
    B_{N,4})}_{L_x^1 L_T^2}\right)^2\right)^{1/2},\label{eq25}
\end{align}
where $B_{N,4}=\abs{\pn}^{k-2}\pn[]\pn\overline{\pn_x}$. By the
H\"older inequality, we have
\begin{align*}
    \norm{B_{N,4}}_{L_x^1 L_T^2}\lesssim
    &\norm{\pn}_{L_x^{k}L_T^\infty}^{k-2}\norm{P_N
    u}_{L_x^{k}L_T^\infty}\norm{\pn\overline{\pn_x}}_{L_x^k L_T^2}\\
    \lesssim &\left[T^{\frac{1}{2}} \norm{u}_{X_T}^k+(1+T^{\frac{1}{4}}\norm{u}_{X_T})\norm{u}_{X_T}^{k-1}\norm{P_{\gg
    1}u}_{X_T}\right]\norm{P_N
    u}_{L_x^{k}L_T^\infty}.
\end{align*}
Thus, the first term in \eqref{eq25} can be bounded by
\begin{align*}
    \lesssim T^{\frac{1}{2}} \norm{u}_{X_T}^k\norm{P_{\gg
    1}u}_{X_T}+(1+T^{\frac{1}{4}}\norm{u}_{X_T})\norm{u}_{X_T}^{k-1}\norm{P_{\gg
    1}u}_{X_T}^2.
\end{align*}

By the H\"older inequality, we get
\begin{align*}
    \norm{B_{N,4}}_{L_x^{\frac{1}{1-\eps}} L_T^2}\lesssim
    &\norm{\pn}_{L_x^{\frac{k(k-2)}{k(1-\eps)-2}}L_T^\infty}^{k-2}\norm{P_N
    u}_{L_x^{k}L_T^\infty}\norm{\pn\overline{\pn_x}}_{L_x^k L_T^2}\\
    \lesssim &\left[T^{\frac{1}{2}} \norm{u}_{X_T}^k+(1+T^{\frac{1}{4}}\norm{u}_{X_T})\norm{u}_{X_T}^{k-1}\norm{P_{\gg
    1}u}_{X_T}\right]\norm{P_N
    u}_{L_x^{k}L_T^\infty}.
\end{align*}
Noticing that $B_{N,4}$ has Fourier support in $|\xi|\sim N$, we can
repeat the procedure which we use to deal with the second term in
\eqref{eq15}, and obtain that the second term in \eqref{eq25} can be
bounded by
\begin{align*}
    \lesssim T^{\frac{1}{2}} \norm{u}_{X_T}^{2k}\norm{P_{\gg
    1}u}_{X_T}+(1+T^{\frac{1}{4}}
   \norm{u}_{X_T})\norm{u}_{X_T}^{2k-1}\norm{P_{\gg
    1}u}_{X_T}^2.
\end{align*}
Therefore, we obtain
\begin{align*}
    &\left(\sum_{N\gg 1}\norm{\int_0^t S(t-\tau)I_{N,4}(\tau)
    d\tau}_{Y_T}^2\right)^{1/2}\\
    \lesssim &(1+\norm{u}_{X_T}^k)\left[T^{\frac{1}{2}} \norm{u}_{X_T}^k\norm{P_{\gg
    1}u}_{X_T}+(1+T^{\frac{1}{4}}\norm{u}_{X_T})\norm{u}_{X_T}^{k-1}\norm{P_{\gg
    1}u}_{X_T}^2\right].
\end{align*}

\subsection{Nonlinear estimates of $I_{N,5}$}

From \eqref{eq5}, \eqref{eq12}, \eqref{eq7} and \eqref{eq8}, we have
\begin{align}
    &\left(\sum_{N\gg 1}\norm{\int_0^t S(t-\tau)\ep B_{N,5}
    d\tau}_{Y_T}^2\right)^{1/2}\no\\
    \lesssim &\left(\sum_{N\gg 1}\norm{ B_{N,5}
    }_{L_x^1 L_T^2}^2\right)^{1/2}\no\\
    &+\left(\sum_{N\gg 1}\left(\sum_{M}\norm{P_M(\ep
    B_{N,5})}_{L_x^1 L_T^2}\right)^2\right)^{1/2},\label{eq26}
\end{align}
where $B_{N,5}=\abs{\pn}^{2k}\pn[]$. By the H\"older inequality, we
have
\begin{align*}
    \norm{B_{N,5}}_{L_x^1 L_T^2}\lesssim
    &T^{\frac{1}{2}}\norm{\pn}_{L_x^{4k}L_T^\infty}^{2k}\norm{P_N
    u}_{L_x^2L_T^\infty}\\
    \lesssim &T^{\frac{1}{2}}\norm{u}_{X_T}^{2k}\norm{P_N
    u}_{L_x^2L_T^\infty},
\end{align*}
and
\begin{align*}
    \norm{B_{N,5}}_{L_x^{\frac{1}{1-\eps}} L_T^2}\lesssim
    &T^{\frac{1}{2}}\norm{\pn}_{L_x^{\frac{4k}{1-2\eps}}L_T^\infty}^{2k}\norm{P_N
    u}_{L_x^2L_T^\infty}\\
    \lesssim &T^{\frac{1}{2}}\norm{u}_{X_T}^{2k}\norm{P_N
    u}_{L_x^2L_T^\infty}.
\end{align*}
Thus, in a similar way as dealing with $I_{N,1}$ and $I_{N,4}$, and
noticing that $B_{N,5}$ has Fourier support in $|\xi|\sim N$, we can
bound \eqref{eq26} by
\begin{align*}
    \lesssim T^{\frac{1}{2}}\norm{u}_{X_T}^{3k}\norm{P_{\gg
    1}u}_{X_T}.
\end{align*}

\section{A priori estimates for solutions}

By the scaling argument
\begin{align*}
    u(t,x)\mapsto
    u_\gamma(t,x)=\frac{1}{\gamma^{1/k}}u(\frac{t}{\gamma^2},\frac{x}{\gamma}),
\end{align*}
we have
\begin{align*}
    \norm{u_{0,\gamma}}_{L^2}=&\gamma^{\frac{1}{2}-\frac{1}{k}}\norm{u_0}_{L^2},\\
    \norm{u_{0,\gamma}}_{\dot{H}^{\frac{1}{2}}}=&\frac{1}{\gamma^{1/k}}\norm{u_0}_{\dot{H}^{\frac{1}{2}}}.
\end{align*}
Thus, we may rescale
\begin{align*}
    \norm{P_{\lesssim 1}u_{0,\gamma}}_{L^2}\ls&
    \gamma^{\frac{1}{2}-\frac{1}{k}}\norm{u_0}_{L^2}=C_{low},\\
    \norm{P_{\gg
    1}u_{0,\gamma}}_{H^{\frac{1}{2}}}\ls&\frac{1}{\gamma^{1/k}}\norm{u_0}_{H^{\frac{1}{2}}}<C_{high}\ll
    1,
\end{align*}
where we choose $\gamma=\gamma(\norm{u_0}_{H^{1/2}})\gg 1$, and take
the time interval $T$ depending on $\gamma$ later. We now drop the
writing of the scaling parameter $\gamma$ and assume
\begin{align*}
    \norm{P_{\lesssim 1}u_0}_{L^2}\ls& C_{low},\\
    \norm{P_{\gg 1}u_0}_{H^{\frac{1}{2}}}\ls& C_{high}\ll 1.
\end{align*}
We now apply this to the norms $X_T$ and $H^{1/2}$, and define new
version of the norms of $X_T$ and $H^{1/2}$, given by with the
decomposition $I=P_{\lesssim 1}+P_{\gg 1}$,
\begin{align*}
    \norm{u}_{\tilde{X}_T}=\frac{1}{C_{low}}\norm{P_{\lesssim
    1}u}_{X_T}+\frac{1}{C_{high}}\norm{P_{\gg 1}u}_{X_T},
\end{align*}
and
\begin{align*}
    \norm{\phi}_{\tilde{H}^{1/2}}=\frac{1}{C_{low}}\norm{P_{\lesssim
    1}\phi}_{L^2}+\frac{1}{C_{high}}\norm{P_{\gg 1}\phi}_{H^{1/2}},
\end{align*}
which implies that $\norm{u_0}_{\tilde{H}^{1/2}}\ls 2$.

For the low frequency part, we have the following estimates.

\begin{lemma}\label{lem6}
Let $u$ be a solution of \eqref{dnls}-\eqref{data}. Then
\begin{align*}
    \norm{P_{\lesssim 1}u}_{X_T}\lesssim
    C_{low}+T^{1/2}\norm{u}_{X_T}^{k+1}.
\end{align*}
\end{lemma}

\begin{proof}
Using the integral equation of \eqref{dnls}
\begin{align*}
    u(t)=S(t)u_0-\lambda\int_0^t S(t-\tau)|u(\tau)|^ku_x(\tau)d\tau,
\end{align*}
and by \eqref{eq1}, \eqref{eq2}, \eqref{eq3}, \eqref{eq9},
\eqref{eq10}, \eqref{eq11} and the H\"older inequality, we have
\begin{align*}
    \norm{P_{\lesssim 1}u}_{X_T}\lesssim & \norm{S(t)P_{\lesssim
    1}u_0}_{X_T}+\norm{\int_0^t S(t-\tau)P_{\lesssim
    1}(|u|^ku_x)(\tau)d\tau}_{X_T}\\
    \lesssim& \norm{P_{\lesssim 1}u_0}_{L^2}+\norm{P_{\lesssim
    1}(|u|^ku_x)}_{L_T^1 H_x^{1/2}}
    \lesssim C_{low}+\norm{|u|^ku_x}_{L_T^1 L_x^2}\\
    \lesssim& C_{low}+T^{1/2}\norm{u}_{L_x^{2k}L_T^\infty}^k\norm{u_x}_{ L_x^\infty
    L_T^2}\\
    \lesssim & C_{low}+T^{1/2}\norm{u}_{X_T}^{k+1},
\end{align*}
which is the desired result.
\end{proof}

For the high frequency part, we have

\begin{lemma}\label{lem7}
Let $u$ and $v_N$ be given in \eqref{transf}. Then
\begin{align*}
    \norm{P_{\gg 1}u}_{X_T}\lesssim (1+\norm{u}_{L_T^\infty
    H_x^{1/2}}^{2k})\left(\sum_{N\gg 1}\norm{v_N}_{Y_T}^2\right)^{1/2}.
\end{align*}
\end{lemma}

\begin{proof}
By \eqref{transf}, we have
\begin{align*}
    P_N u=e^{\frac{i\lambda}{2}\int_{-\infty}^x |P_{\ll N}u|^k
    dy}v_N.
\end{align*}
For $L_T^\infty H_x^{1/2}$-norm, by the interpolation theorem, we
obtain for $N\gg 1$,
\begin{align*}
    \norm{P_N u}_{H_x^{1/2}}\lesssim & \norm{P_N u}_{L^2}^{\frac{1}{2}}\norm{P_N u}_{H^1}^{\frac{1}{2}}
    \lesssim \norm{v_N}_{L^2}^{\frac{1}{2}}\left(\norm{P_N u}_{L^2}+\norm{\partial_xP_N u}_{L^2}\right)^{\frac{1}{2}}\\
    \lesssim &\norm{v_N}_{L^2}^{\frac{1}{2}}\left(\norm{v_N}_{L^2}
    +\norm{\abs{P_{\ll N}u}^k v_N}_{L^2}+\norm{\partial_xv_N}_{L^2}\right)^{\frac{1}{2}}\\
    \lesssim & \norm{v_N}_{L^2}^{\frac{1}{2}}\left(\norm{P_{\ll N}u}_{L_x^{4k}}^k
    \norm{v_N}_{L_x^4}+\norm{v_N}_{H^1}\right)^{\frac{1}{2}}\\
    \lesssim & \left(1+\norm{P_{\ll
    N}u}_{H_x^{1/2}}^k\right)^{\frac{1}{2}}\norm{v_N}_{H_x^{1/2}}
    \lesssim  \left(1+\norm{P_{\ll
    N}u}_{H_x^{1/2}}^k\right)\norm{v_N}_{H_x^{1/2}},
\end{align*}
which yields the desired estimate by summing on $l_N^2$.

For the $L_x^\infty L_T^2$-norm, noticing that
\begin{align*}
    \partial_x P_Nu=e^{\frac{i\lambda}{2}\int_{-\infty}^x |P_{\ll N}u|^k
    dy}(\partial_x v_N+\frac{i\lambda}{2}|\pn|^k v_N),
\end{align*}
we have
\begin{align}
    &\norm{\partial_x P_Nu}_{L_x^\infty L_T^2}\lesssim
    \norm{\partial_x v_N}_{L_x^\infty L_T^2}\no\\
    &\quad +\norm{\tilde{P}_N\left(\sum_{N_1}P_{N_1}(e^{\frac{i\lambda}{2}\int_{-\infty}^x |P_{\ll N}u|^k
    dy})|\pn|^k\sum_{N_2}P_{N_2}v_N \right)}_{L_x^\infty
    L_T^2}.\label{eq27}
\end{align}
To estimate the second term \eqref{eq27}, we split the sum
$\sum_{N_2}=\sum_{N_2\sim N}+\sum_{N_2\nsim N}$. For $N_2\sim N$,
from the Bernstein inequality, we bound \eqref{eq27} by
\begin{align*}
    \lesssim & \norm{|\pn|^k\sum_{N_2\sim N}P_{N_2}v_N }_{L_x^\infty
    L_T^2}\lesssim \norm{\pn}_{L_{x,T}^\infty}^k\sum_{N_2\sim
    N}\norm{P_{N_2}v_N}_{L_x^\infty L_T^2}\\
    \lesssim &N\norm{D_x^{-1/k}\pn}_{L_{x,T}^\infty}^k\sum_{N_2\sim
    N}\norm{ P_{N_2}v_N}_{L_x^\infty L_T^2}\\
    \lesssim &\norm{\pn}_{L_T^\infty H_x^{1/2}}^k\sum_{N_2\sim
    N}\norm{P_{N_2}\partial_xv_N}_{L_x^\infty L_T^2}\\
    \lesssim &\norm{u}_{L_T^\infty H_x^{1/2}}^k\sum_{N_2\sim
    N}\norm{P_{N_2}\partial_xv_N}_{L_x^\infty L_T^2}.
\end{align*}
For the part $N_2\nsim N$, we split it as $\sum_{N_2\nsim
N}=\sum_{N_2\ll N}+\sum_{N_2\gg N}$. Noticing that for $N_2\ll N$,
\begin{align}\label{eq28}
    &\tilde{P}_N\left(\sum_{N_1} P_{N_1}(e^{\frac{i\lambda}{2}\int_{-\infty}^x |P_{\ll N}u|^k
    dy})|\pn|^k\sum_{N_2}P_{N_2}v_N \right)\\
    =&\tilde{P}_N\left(P_{\sim N}(e^{\frac{i\lambda}{2}\int_{-\infty}^x |P_{\ll N}u|^k
    dy})|\pn|^kP_{\ll N}v_N \right),\no
\end{align}
and for $N_2\gg N$,
\begin{align*}
    \eqref{eq28}=\tilde{P}_N\left(\sum_{N_1\sim N_2\gg N}P_{N_1}(e^{\frac{i\lambda}{2}\int_{-\infty}^x |P_{\ll N}u|^k
    dy})|\pn|^kP_{N_2}v_N \right),
\end{align*}
we can bound \eqref{eq27}, in view of the Bernstein inequality and
the H\"older inequality, by
\begin{align*}
    \lesssim & \norm{P_{\sim N}e^{\frac{i\lambda}{2}\int_{-\infty}^x |P_{\ll N}u|^k
    dy}}_{L_{x,T}^\infty}\norm{\pn}_{L_{x,T}^\infty}^k\norm{P_{\ll
    N}v_N}_{L_{x,T}^\infty}\\
    &+\sum_{N_1\sim N_2\gg N}\norm{P_{N_1}e^{\frac{i\lambda}{2}\int_{-\infty}^x |P_{\ll N}u|^k
    dy}}_{L_{x,T}^\infty}\norm{\pn}_{L_{x,T}^\infty}^k\norm{P_{N_2}v_N}_{L_{x,T}^\infty}\\
    \lesssim &N^{-\frac{k}{(2k+1)}}\norm{P_{\sim N}\partial_x e^{\frac{i\lambda}{2}\int_{-\infty}^x |P_{\ll N}u|^k
    dy}}_{L_{x,T}^\infty}\\
    &\qquad\cdot\norm{D_x^{-\frac{1}{(2k+1)}}\pn}_{L_{x,T}^\infty}^k\norm{D_x^{-\frac{1}{(2k+1)}}P_{\ll
    N}v_N}_{L_{x,T}^\infty}\\
    &+\!\!\sum_{N_1\sim N_2\gg N}\!\! N_1^{-\frac{k}{(2k+1)}}
    \norm{P_{N_1}\partial_x e^{\frac{i\lambda}{2}\int_{-\infty}^x |P_{\ll N}u|^k
    dy}}_{L_{x,T}^\infty}\\
    &\qquad\cdot\norm{D_x^{-\frac{1}{(2k+1)}}\pn}_{L_{x,T}^\infty}^k
    \norm{D_x^{-\frac{1}{(2k+1)}}P_{N_2}v_N}_{L_{x,T}^\infty}\\
    \lesssim
    &N^{-\frac{k}{(2k+1)}}\norm{\pn}_{L_{x,T}^\infty}^k\norm{u}_{L_T^\infty
    H_x^{1/2}}^k\norm{v_N}_{L_T^\infty H_x^{1/2}}\\
    &+\!\!\! \sum_{N_1\sim N_2\gg N}\!\!\! N_1^{-\frac{k}{3(2k+1)}}N_2^{-\frac{k}{3(2k+1)}}N^{-\frac{k}{3(2k+1)}}
    \norm{\pn}_{L_{x,T}^\infty}^k\norm{u}_{L_T^\infty
    H_x^{1/2}}^k\norm{P_{N_2}v_N}_{L_T^\infty H_x^{1/2}}\\
    \lesssim &\norm{u}_{L_T^\infty
    H_x^{1/2}}^{2k}\norm{v_N}_{L_T^\infty H_x^{1/2}}.
\end{align*}
Therefore, summing on $l_N^2$, we complete the proof for the
$L_x^\infty L_T^2$-norm.

For the $L_x^2L_T^\infty$-norm, it is easy to obtain the desired
result since $|P_Nu|=|v_N|$.

We turn to estimate the $L_x^4 L_T^\infty$-norm. It is similar with
the proof for the $L_x^\infty L_T^2$-norm, since
$\norm{\dx^{1/4}P_Nu}_{L_x^4 L_T^\infty}\sim
N^{1/4}\norm{P_Nu}_{L_x^4 L_T^\infty}$ for $N\gg 1$. In fact, we
have
\begin{align}\label{eq29}
    \norm{\dx^{1/4}P_Nu}_{L_x^4 L_T^\infty}\sim N^{1/4}
    \norm{\sum_{N_1}P_{N_1}(e^{\frac{i\lambda}{2}\int_{-\infty}^x |P_{\ll N}u|^k
    dy})\sum_{N_2}P_{N_2}v_N}_{L_x^4 L_T^\infty}.
\end{align}
We also split $\sum_{N_2}=\sum_{N_2\sim N}+\sum_{N_2\nsim N}$. For
$N_2\sim N$, we bound \eqref{eq29} by
\begin{align*}
    \lesssim &N^{1/4}\sum_{N_2\sim N}\norm{P_{N_2}v_N}_{L_x^4
    L_T^\infty}\lesssim \norm{\dx^{1/4}v_N}_{L_x^4 L_T^\infty}.
\end{align*}
For the part $N_2\nsim N$, we split it as $\sum_{N_2\nsim
N}=\sum_{N_2\ll N}+\sum_{N_2\gg N}$. Noticing that for $N_2\ll N$,
\begin{align}\label{eq30}
    &\tilde{P}_N\left(\sum_{N_1} P_{N_1}(e^{\frac{i\lambda}{2}\int_{-\infty}^x |P_{\ll N}u|^k
    dy})\sum_{N_2}P_{N_2}v_N \right)\\
    =&\tilde{P}_N\left(P_{\sim N}(e^{\frac{i\lambda}{2}\int_{-\infty}^x |P_{\ll N}u|^k
    dy})P_{\ll N}v_N \right),\no
\end{align}
and for $N_2\gg N$,
\begin{align*}
    \eqref{eq30}=\tilde{P}_N\left(\sum_{N_1\sim N_2\gg N}P_{N_1}(e^{\frac{i\lambda}{2}\int_{-\infty}^x |P_{\ll N}u|^k
    dy})P_{N_2}v_N \right),
\end{align*}
we can bound \eqref{eq29}, in view of the Bernstein inequality and
the H\"older inequality, by
\begin{align*}
    \lesssim &\norm{\dx^{1/4}v_N}_{L_x^4 L_T^\infty}
    +N^{1/4}\norm{P_{\sim N}(e^{\frac{i\lambda}{2}\int_{-\infty}^x |P_{\ll N}u|^k
    dy})}_{L_{x,T}^\infty}\norm{P_{\ll N}v_N}_{L_x^4L_T^\infty}\\
    &+N^{1/4}\sum_{N_1\sim N_2\gg N}\norm{P_{N_1}(e^{\frac{i\lambda}{2}\int_{-\infty}^x |P_{\ll N}u|^k
    dy})}_{L_{x,T}^\infty}\norm{P_{N_2}v_N}_{L_x^4L_T^\infty}\\
    \lesssim &\norm{\dx^{1/4}v_N}_{L_x^4
    L_T^\infty}+N^{-3/4}\norm{\pn}_{L_{x,T}^\infty}^k\norm{\dx^{1/4}P_{\ll
    N}v_N}_{L_x^4L_T^\infty}\\
    &+\sum_{N_1\sim N_2\gg N}N_1^{-1/3}N_2^{-1/3}N^{-1/3} \norm{\pn}_{L_{x,T}^\infty}^k
    \norm{\dx^{1/4}P_{N_2}v_N}_{L_x^4L_T^\infty}\\
    \lesssim &(1+\norm{u}_{L_T^\infty
    H_x^{1/2}}^k)\norm{\dx^{1/4}v_N}_{L_x^4L_T^\infty},
\end{align*}
which yields the desired estimate by applying $l_N^2$-sum.

Thus, we complete the proof of this Lemma.
\end{proof}

Of course, we need the following estimate of the data.

\begin{lemma}\label{lem8}
For any $u_0\in H^{1/2}$, we have
\begin{align}\label{eq31}
    \left(\sum_{N\gg 1}\norm{S(t)(e^{-\frac{i\lambda}{2}\int_{-\infty}^x|P_{\ll N}u_0|^k dy}
    P_Nu_0)}_{Y_T}^2\right)^{1/2}\lesssim
    (1+\norm{u_0}_{H^{1/2}}^k)\norm{P_{\gg 1}u_0}_{H^{1/2}}.
\end{align}
\end{lemma}

\begin{proof}
From \eqref{eq1}, \eqref{eq2} and \eqref{eq3}, we bound the left
hand side of \eqref{eq31} by
\begin{align}\label{eq32}
    \lesssim &\left(\sum_{N\gg 1}\norm{e^{-\frac{i\lambda}{2}\int_{-\infty}^x|P_{\ll N}u_0|^k dy}
    P_Nu_0}_{H^{1/2}}^2\right)^{1/2}\\
    &+\left(\sum_{N\gg
    1}\left(\sum_M\norm{P_M(e^{-\frac{i\lambda}{2}\int_{-\infty}^x|P_{\ll N}u_0|^k dy}
    P_Nu_0)}_{H^{1/2}}\right)^2\right)^{1/2}. \label{eq33}
\end{align}

From Lemma~\ref{lem1}, we have
\begin{align*}
    \eqref{eq32}\lesssim &\left(\sum_{N\gg 1}[\norm{P_{\ll N}u_0}_{L^{4k}}^{2k}
    \norm{P_Nu_0}_{L^4}^2+\norm{P_Nu_0}_{H^{1/2}}^2]\right)^{1/2}\\
    \lesssim &(1+\norm{u_0}_{H^{1/2}}^2)\norm{P_{\gg
    1}u_0}_{H^{1/2}}.
\end{align*}
For the second term \eqref{eq33}, it is similar with \eqref{eq17}.
We split the sum $\sum_M=\sum_{M\lesssim N}+\sum_{M\gg N}$. By the
Bernstein inequality, the H\"older inequality and the Sobolev
embedding theorem, we bound \eqref{eq33} by
\begin{align*}
    \lesssim &\left(\sum_{N\gg 1}\left(\sum_{M\lesssim N}\langle
    M\rangle^{\frac{1}{2}}\norm{P_Nu_0}_{L_x^2}\right)^2\right)^{1/2}\\
    &+\left(\sum_{N\gg
    1}\left(\sum_{M\gg N}\norm{\dx^{\frac{1}{2}}P_M(e^{-\frac{i\lambda}{2}\int_{-\infty}^x|P_{\ll N}u_0|^k dy}
    P_Nu_0)}_{L_x^2}\right)^2\right)^{1/2}\\
    \lesssim &\left(\sum_{N\gg
    1}\left(N^{1/2}\norm{P_Nu_0}_{L_x^2}\right)^2\right)^{1/2}\\
    &+\left(\sum_{N\gg
    1}\left(\sum_{M\gg N}M^{\frac{1}{2}}\norm{(P_{\sim M}e^{-\frac{i\lambda}{2}\int_{-\infty}^x|P_{\ll N}u_0|^k
    dy})
    P_Nu_0}_{L_x^2}\right)^2\right)^{1/2}\\
    \lesssim &\norm{P_{\gg 1}u_0}_{H^{1/2}}\\
    &+\left(\sum_{N\gg
    1}\left(\sum_{M\gg N}M^{\frac{1}{2}}\norm{P_{\sim M}e^{-\frac{i\lambda}{2}\int_{-\infty}^x|P_{\ll N}u_0|^k
    dy}}_{L^4}\norm{P_Nu_0}_{L_x^4}\right)^2\right)^{1/2}\\
    \lesssim &\norm{P_{\gg 1}u_0}_{H^{1/2}}\\
    &+\left(\sum_{N\gg
    1}\left(\sum_{M\gg N}M^{-\frac{1}{2}}\norm{P_{\sim M}\partial_xe^{-\frac{i\lambda}{2}\int_{-\infty}^x|P_{\ll N}u_0|^k
    dy}}_{L^4}\norm{P_Nu_0}_{L_x^4}\right)^2\right)^{1/2}\\
    \lesssim &\norm{P_{\gg 1}u_0}_{H^{1/2}}+\left(\sum_{N\gg
    1}\left(\norm{P_{\ll N}u_0}_{L^{4k}}^k
    \norm{P_Nu_0}_{L_x^4}\right)^2\right)^{1/2}\\
    \lesssim &(1+\norm{u_0}_{H^{1/2}}^k)\norm{P_{\gg
    1}u_0}_{H^{1/2}},
\end{align*}
which yields the desired result.
\end{proof}

With the help of the above lemmas, we can prove the following
proposition which yields the a priori estimate.

\begin{proposition}\label{prop3}
Let $u$ be a smooth solution to \eqref{dnls}-\eqref{data} and
$0<T\ls C_{high}^4$. Then we have
\begin{align*}
    \norm{u}_{\tilde{X}_T}\ls
    C(C_{low})+C(C_{low}+\norm{u}_{\tilde{X}_T})^{3k}(T^{1/4}
    +C_{high})\norm{u}_{\tilde{X}_T}.
\end{align*}
\end{proposition}

\begin{proof}
Noticing that
\begin{align*}
    \norm{P_{\lesssim 1}u}_{X_T}\lesssim
    C_{low}\norm{u}_{\tilde{X}_T}, \quad \norm{P_{\gg 1}u}_{X_T}\lesssim
    C_{high}\norm{u}_{\tilde{X}_T}.
\end{align*}
and from Lemmas \ref{lem6}, \ref{lem7}, \ref{lem8} and Proposition
\ref{prop2}, we obtain through a complicated computation
\begin{align*}
&\norm{u}_{\tilde{X}_T}=\frac{1}{C_{low}}\norm{P_{\lesssim
    1}u}_{X_T}+\frac{1}{C_{high}}\norm{P_{\gg 1}u}_{X_T}\\
    \lesssim &1+\frac{1}{C_{low}}T^{\frac{1}{2}}\norm{u}_{X_T}^{k+1}
    +\frac{1}{C_{high}}(1+\norm{u}_{L_T^\infty
    H_x^{1/2}}^{2k})\left(\sum_{N\gg
    1}\norm{v_N}_{Y_T}^2\right)^{1/2}\\
    \lesssim &1+\frac{1}{C_{low}}T^{\frac{1}{2}}\norm{u}_{X_T}^{k+1}
    +\frac{1}{C_{high}}(1+\norm{u}_{L_T^\infty
    H_x^{1/2}}^{2k})(1+\norm{u_0}_{H^{1/2}}^k)\norm{P_{\gg
    1}u_0}_{H^{1/2}}\\
    &+\frac{1}{C_{high}}(1+\norm{u}_{L_T^\infty
    H_x^{1/2}}^{2k})\Big\{T^{\frac{1}{2}}\norm{u}_{X_T}^{3k}\norm{P_{\gg 1}
    u}_{X_T}\\
    &+(1+\norm{u}_{X_T}^k)\left[T^{\frac{1}{2}}\norm{u}_{X_T}^{\tilde{k}+1}\norm{P_{\gg
1}u}_{X_T}^{k-\tilde{k}}
    +(1+T^{\frac{1}{4}}\norm{u}_{X_T})\norm{u}_{X_T}^{\tilde{k}}\norm{P_{\gg
    1}u}_{X_T}^{k+1-\tilde{k}}\right]\\
    &+T^{\frac{1}{2}}\left(\norm{u}_{X_T}^{2k-1}
    +\norm{u}_{X_T}^{(5k-2)/2}\right)\norm{P_{\gg
    1}u}_{X_T}\\
    &+T^{\frac{1}{4}}(1+T^{\frac{1}{4}}\norm{u}_{X_T})^{\frac{1}{2}}\norm{u}_{X_T}^{2k-1}\norm{P_{\gg
    1}u}_{X_T}^{\frac{3}{2}}\\
    &+(1+\norm{u}_{X_T}^k)\left[T^{\frac{1}{2}} \norm{u}_{X_T}^k\norm{P_{\gg
    1}u}_{X_T}+(1+T^{\frac{1}{4}}\norm{u}_{X_T})^2\norm{u}_{X_T}^{k-1}\norm{P_{\gg
    1}u}_{X_T}^2\right]\Big\}\\
    \ls &C+C(C_{low})T^{\frac{1}{2}}\norm{u}_{\tilde{X}_T}^{k+1}
    +(1+\norm{u}_{L_T^\infty
    H_x^{1/2}}^{2k})(1+\norm{u_0}_{H^{1/2}}^k)\\
    &+C(C_{low})(1+\norm{u}_{L_T^\infty
    H_x^{1/2}}^{2k})\Big\{T^{\frac{1}{4}}\norm{u}_{\tilde{X}_T}^{3k+1}+T^{\frac{1}{4}}
    (1+\norm{u}_{\tilde{X}_T}^k)\norm{u}_{\tilde{X}_T}^{k+1}\\
    &+T^{\frac{1}{2}}(\norm{u}_{\tilde{X}_T}^{k+1}
    +\norm{u}_{\tilde{X}_T}^{5k/2})
    +C_{high}^{k-\tilde{k}}(1+T^{\frac{1}{4}}\norm{u}_{\tilde{X}_T})(1+\norm{u}_{\tilde{X}_T}^k)
    \norm{u}_{\tilde{X}_T}^{k+1}\\
    &+C_{high}^{3/2}(1+T^{\frac{1}{4}}\norm{u}_{\tilde{X}_T})\norm{u}_{\tilde{X}_T}^{2k}
    +C_{high}^{2}(1+T^{\frac{1}{4}}\norm{u}_{\tilde{X}_T})^2(1+\norm{u}_{\tilde{X}_T}^k)
    \norm{u}_{\tilde{X}_T}^{k+1}\Big\}.
\end{align*}
Notice that
\begin{align*}
    \norm{u(t)}_{H^{1/2}}\lesssim \norm{P_{\lesssim
    1}u(t)}_{L^2}+C_{high}\norm{P_{\gg 1}}_{\tilde{H}^{1/2}}.
\end{align*}
The high frequency part $C_{high}\norm{P_{\gg 1}}_{\tilde{H}^{1/2}}$
can be absorbed into the $\tilde{X}_T$-norm. Then substituting Lemma
\ref{lem6} again in estimating the low frequency part of the norm
$\norm{P_{\lesssim 1}u}_{L_T^\infty H_x^{1/2}}$, we complete the
proof of Proposition \ref{prop3}.
\end{proof}

From Proposition \ref{prop3}, we have the following a priori
estimate for the solution of \eqref{dnls}-\eqref{data} if we take
$T$ and $C_{high}$ small enough.

\begin{corollary}
Let $u$ be a smooth solution to \eqref{dnls}-\eqref{data}. we have
\begin{align*}
    \norm{u}_{\tilde{X}_T}\lesssim C_{low}+C_{high},
\end{align*}
for $T$ and $C_{high}$ small enough.
\end{corollary}

For the proof of Theorem~\ref{thm}, we can follow the compactness
argument with the a priori estimate. Since the proof is standard, we
omit the details and refer to the papers
\cite{KT05,KT03,MR04a,MR04b,Po91,Tao04}.

\section*{Acknowledgements} The author would like to thank the referees for valuable
comments and suggestions on the original manuscript and Prof. L.
Hsiao for her frequent encouragement.

\medskip

Received July 25, 2006; revised November 2006.

\medskip

\end{document}